%% file: main.tex
\documentclass[final,5p,times,twocolumn]{elsarticle}

\usepackage[T1]{fontenc}
\usepackage[utf8]{inputenc}
\usepackage{lmodern}
\usepackage{microtype}
\usepackage{amsmath,amssymb,amsfonts,mathtools,bm}
\usepackage{booktabs,tabularx,array}
\usepackage{xcolor}
\usepackage{graphicx}
\usepackage{subcaption}
\usepackage[hidelinks]{hyperref}
\usepackage{cleveref}
\usepackage{enumitem}
\usepackage{placeins}
\input{qpde_figure_style.tex}

\journal{Computer Physics Communications}

\newcommand{\bA}{\bm{A}}
\newcommand{\bD}{\bm{D}}
\newcommand{\bE}{\bm{E}}
\newcommand{\bF}{\bm{F}}
\newcommand{\bH}{\bm{H}}
\newcommand{\bI}{\bm{I}}
\newcommand{\bL}{\bm{L}}
\newcommand{\bP}{\bm{P}}
\newcommand{\bQ}{\bm{Q}}
\newcommand{\bS}{\bm{S}}
\newcommand{\bT}{\bm{T}}
\newcommand{\bb}{\bm{b}}
\newcommand{\be}{\bm{e}}
\newcommand{\bq}{\bm{q}}
\newcommand{\bu}{\bm{u}}
\newcommand{\bx}{\bm{x}}
\newcommand{\by}{\bm{y}}
\newcommand{\bY}{\bm{Y}}
\newcommand{\R}{\mathbb{R}}

\newcommand{\rev}[1]{\textcolor{black}{#1}}

\begin{document}

\begin{frontmatter}

\title{\rev{Reusable QUBO-Derived Block-Schur Inverse Factors for Fixed-Operator Finite-Difference PDE Solvers}}

\author[aff1]{Arash Fath Lipaei\corref{cor1}}
\ead{arash.fa.li.2000@gmail.com}
\cortext[cor1]{Corresponding author.}

\author[aff2]{\"Ozg\"ur E. M\"ustecapl{\i}o\u{g}lu}
\author[aff3,aff4]{\"Ozlem Salehi}
\author[aff5]{Selin Aslan}

\affiliation[aff1]{organization={Computational Sciences and Engineering Department, Ko\c{c} University},
                   addressline={Rumelifeneri Yolu},
                   city={Istanbul},
                   postcode={34450},
                   country={T\"urkiye}}

\affiliation[aff2]{organization={Department of Physics, Ko\c{c} University},
                   addressline={Rumelifeneri Yolu},
                   city={Istanbul},
                   postcode={34450},
                   country={T\"urkiye}}
\affiliation[aff3]{organization={Algorithmiq S.r.l.},
                   city={Milano},
                   postcode={20123},
                   country={Italy}}
\affiliation[aff4]{organization={Faculty of Engineering and Natural Sciences, Sabancı University},
                    addressline={Üniversite Caddesi},
                   city={Istanbul},
                   postcode={34956},
                   country={T\"urkiye}}
\affiliation[aff5]{organization={Department of Mathematics, Ko\c{c} University},
                   addressline={Rumelifeneri Yolu},
                   city={Istanbul},
                   postcode={34450},
                   country={T\"urkiye}}
\begin{abstract}
Many simulations of partial differential equations require the repeated solution of linear systems whose matrix remains fixed while the right-hand side changes. We present a hybrid offline--online framework for this setting. The discretized operator is organized into small coupled blocks and reduced through block Schur complements. Instead of expressing the complete PDE solution as one large binary optimization problem, quadratic unconstrained binary optimization (QUBO) is used only to construct reusable local inverse factors. These factors are computed during an offline stage, decoded and stored, and then reused in an ordinary classical block solve for subsequent right-hand sides. The optimizer-facing problem size is governed by the local block size and chosen numerical precision, rather than by the dimension of the full spatial grid. We develop the mathematical formulation, analyze the accuracy and resource requirements of the binary representation, and examine the method on representative finite-difference PDE problems, including execution of the local optimization stage on quantum optimization hardware. The resulting architecture provides a reusable interface between conventional numerical linear algebra and QUBO-based quantum optimization.
\end{abstract}

\begin{keyword}
partial differential equations \sep finite differences \sep QUBO \sep Schur complement \sep hybrid computing \sep computational physics \sep quantum optimization
\end{keyword}

\end{frontmatter}

\section{Introduction and related work}
\label{sec:introduction}

Numerical models of diffusion, waves, transport, and many other physical processes reduce, after spatial and temporal discretization, to linear systems of the form
\begin{equation}
    \bA\bx = \bb .
    \label{eq:intro_linear_system}
\end{equation}
In a large class of simulations the operator \(\bA\) remains unchanged over many solves, whereas the right-hand side \(\bb\) changes with time, forcing, boundary data, or an outer iteration. This fixed-operator structure appears naturally in implicit and semi-implicit time integration, multiple-source calculations, parameter studies, and repeated forward solves used in inverse problems. Classical numerical linear algebra already exploits this reuse effectively through sparse factorizations, block elimination, preconditioners, Krylov methods, and multigrid techniques \cite{leveque2007,saad2003,golubvanloan2013,higham2002}. These methods remain the appropriate performance baseline for PDE computation.

We isolate the binary-optimization component of a repeated PDE calculation as a compact, operator-dependent preprocessing task. A direct binary encoding of every degree of freedom makes the optimization problem grow with the full grid, even when the same matrix is reused for hundreds or thousands of right-hand sides. The hybrid formulation preserves the structure exploited by classical solvers and restricts binary optimization to a local task whose cost is distributed over subsequent solves.

Many computational-physics workflows repeatedly solve the same discretized operator for changing source terms, successive time levels at fixed time-step size, or outer iterations. The present construction isolates a bounded local preprocessing kernel for quantum optimization hardware while retaining the global PDE solve in conventional numerical linear algebra. This allows specialized optimization hardware to participate in simulations whose total spatial dimension is much larger than the optimizer-facing problem.

Quadratic unconstrained binary optimization provides a natural interface for such a construction. A QUBO has the form
\begin{equation}
    \min_{\bq\in\{0,1\}^{n}} \; \bq^{\mathsf T}\bQ\bq ,
    \label{eq:intro_qubo}
\end{equation}
where linear terms can be incorporated into the diagonal of $\bQ$, while constant terms can be omitted since they do not affect the optimization procedure~\cite{lucas2014,glover2018tutorial,kochenberger2014}. QUBO and the equivalent Ising formulation are prominent abstractions for quantum annealing~\cite{kadowaki1998,albash2018,dwave_qubo_ising}and have also become widely used across classical, quantum-inspired, and physics-inspired optimization approaches, including digital annealing~\cite{sao2019digital}, simulated bifurcation~\cite{goto2019simulated}, and specialized QUBO solvers such as VeloxQ~\cite{pawlowski2025veloxq}. Alternative quantum optimization architectures, such as the entropy-computing-based Dirac-3~\cite{emami2025entropy}, further broaden this landscape. Real-valued least-squares and linear-algebra problems can be connected to this interface by representing unknown coefficients with binary fixed-point variables and expanding a squared residual into a quadratic polynomial. Previous studies have used this idea for least squares, linear systems, division, and matrix inversion \cite{omalley2016,borle2019,borle2022,jun2021,rogers2020}.

Differential equations have also been approached by reformulating their numerical solutions as QUBO or Ising minimization problems for quantum and quantum-inspired optimization. Existing work includes direct or blockwise solution of discretized systems, annealing formulations of finite-element residuals, coefficient optimization in reduced or basis representations, and quantum or quantum-inspired treatments of fluid and transport problems \cite{srivastava2019,pollachini2021,raisuddin2022feqa,criado2022qade,ray2022fluid,kudo2024,farghadan2024,bazarkhanova2025}. Related work on QUBO sparse approximate inverses further shows that inverse-like operators can themselves be useful optimization targets \cite{suresh2023spai}. More broadly, PDE simulation has also motivated variational, non-conservative-system, and quantum-signal-processing approaches in the circuit model \cite{jaksch2022cfd,sato2024nonconservative,kim2026qsp}. These studies connect differential-equation and linear-algebra calculations to quantum or binary optimization. The present method specifically targets inverse columns of local Schur blocks and reuses the cached factors in a deterministic classical online block solve for changing right-hand sides.

We retain the classical block-Schur structure of the discretized operator and use QUBO only to construct columns of local Schur-block inverses. Once decoded, these inverse columns are assembled into reusable local inverse factors and cached. Every subsequent PDE solve is then a deterministic classical online block solve. The resulting separation is conceptually simple: the offline QUBO preprocessing depends on the fixed operator, while the classical online block solve depends on the changing right-hand side. Consequently, the logical QUBO dimension is controlled by local block size and numerical precision rather than by the total number of unknowns in the PDE grid. The $B=2$ cases are not proposed as competitive alternatives to direct classical inversion as they are prepared for minimum budget quantum optimization experiments.

Logical QUBO size and physical hardware cost are different resource measures. Dense logical interactions may incur substantial overhead when mapped to hardware with restricted connectivity through embedding, routing, chains, or other device-specific transformations~\cite{choi2008,choi2011,cai2014,zbinden2020,konz2021embedding}. Embedding parameters such as chain strength and coefficient rescaling can further affect solution quality~\cite{gilbert2024}. Reducing the logical QUBO dimension therefore has value even when the local QUBO remains dense. We therefore track logical variable count, coupling density, fixed-point range, local inverse residuals, downstream PDE error, and optimizer quality separately.

{\color{black}For fixed-operator PDE problems, the cost of constructing approximate inverse factors can be spread across solves with multiple right-hand sides. Constructing these factors requires balancing the size of the binary representation, its numerical range and precision, and the accuracy of the resulting PDE solution. These considerations motivate a block-Schur construction in which local inverse factors are optimized once and reused in subsequent classical solves.} The contribution of this paper is a reusable block-Schur inverse-factor architecture for fixed-operator finite-difference systems, together with an error and resource analysis that connects fixed-point precision to local inverse accuracy and to the final PDE solution. The numerical study is designed in a way that will analyze different aspects of the method. Exact small-QUBO experiments isolate quantization behavior; multi-PDE tests assess the complete cached solver; larger block studies examine representability and optimization difficulty; and a hardware experiment tests the same local QUBO construction on quantum optimization hardware. Binary optimization enters only in the offline QUBO preprocessing; subsequent solves use the classical online block solve.

Table~\ref{tab:architecture_comparison} summarizes how the proposed architecture differs from direct full-system QUBO formulations and generic QUBO decompositions.

\begin{table*}[t]
\centering
\small
\caption{Comparison of the reusable block-Schur QUBO architecture with direct full-system and generic decomposition strategies. The distinction is the object exposed to optimization and whether operator-dependent work can be reused across changing right-hand sides.}
\label{tab:architecture_comparison}
\begin{tabularx}{\textwidth}{@{}p{3.1cm}XXX@{}}
\toprule
 & Direct full-system QUBO & Generic QUBO decomposition & Reusable block-Schur QUBO \\
\midrule
Object treated & Full PDE solution vector & Large global QUBO & Local Schur inverse columns \\
Dependence on RHS & Changes with $\bb$ & Global instance changes with $\bb$ & Independent of $\bb$ for a fixed operator \\
Reuse & No cached inverse factors & Mainly within one optimization instance & Cached across many RHS/time steps \\
Largest optimizer-facing problem & $N(1+M+K)$ & Set by the decomposition strategy & $B(1+M+K)$ \\
Classical role & Assembly and post-processing & Coordinates decomposed subproblems & Performs deterministic classical online block solve \\
\bottomrule
\end{tabularx}
\end{table*}

Section~\ref{sec:method} develops the block-Schur QUBO formulation and offline--online algorithm; Section~\ref{sec:analysis} derives the error and resource scaling; Section~\ref{sec:experiments} defines the numerical and hardware protocol; Sections~\ref{sec:results}--\ref{sec:conclusion} present the results, limitations, and conclusions.

\section{Reusable block-Schur QUBO formulation}
\label{sec:method}

\subsection{Fixed-operator finite-difference systems}

Consider a PDE whose semi-discrete form can be written as
\begin{equation}
    \frac{d\bu}{dt}=\bL\bu+\bm{N}(\bu,t)+\bm{f}(t),
    \label{eq:semidiscrete_system}
\end{equation}
where \(\bL\) is a linear spatial operator and \(\bm{N}\) collects terms treated explicitly. A first-order implicit--explicit step gives
\begin{equation}
    \underbrace{(\bI-\Delta t\,\bL)}_{\bA}\bu^{n+1}
    =\bu^n+\Delta t\,\bm{N}(\bu^n,t_n)+\Delta t\,\bm{f}^{n+1}
    \equiv \bb^n .
    \label{eq:fixed_operator}
\end{equation}
When the spatial discretization and time step are fixed, \(\bA\) is unchanged and only \(\bb^n\) varies. Steady elliptic and frequency-domain problems lead directly to the same algebraic form \(\bA\bu=\bb\). The construction below therefore depends on a matrix property---a reusable fixed operator with a compatible block ordering---rather than on a particular PDE class.

Partition the \(N=n_bB\) unknowns into \(n_b\) blocks of size \(B\), \(\bu=[\bu_1^{\mathsf T},\ldots,\bu_{n_b}^{\mathsf T}]^{\mathsf T}\). For the finite-difference systems considered here, a natural ordering produces the block-tridiagonal system
\begin{equation}
\begin{bmatrix}
\bD_1 & \bF_1 \\
\bE_2 & \bD_2 & \bF_2 \\
& \ddots & \ddots & \ddots \\
&& \bE_{n_b} & \bD_{n_b}
\end{bmatrix}
\begin{bmatrix}
\bu_1\\ \bu_2\\ \vdots\\ \bu_{n_b}
\end{bmatrix}
=
\begin{bmatrix}
\bb_1\\ \bb_2\\ \vdots\\ \bb_{n_b}
\end{bmatrix}.
\label{eq:blocktri}
\end{equation}
Dirichlet conditions are incorporated in the usual finite-difference assembly.

\subsection{Reusable block-Schur inverse factors}

Applying block Gaussian elimination to Eq.~\eqref{eq:blocktri} gives the Schur recursion
\begin{align}
    \bS_1 &= \bD_1, \label{eq:schur_first}\\
    \bT_1 &= \bS_1^{-1}\bF_1,\\
    \bS_i &= \bD_i-\bE_i\bT_{i-1},
    \qquad i=2,\ldots,n_b-1, \label{eq:schur_middle}\\
    \bT_i &= \bS_i^{-1}\bF_i,
    \qquad i=2,\ldots,n_b-1,\\
    \bS_{n_b} &= \bD_{n_b}-\bE_{n_b}\bT_{n_b-1}. \label{eq:schur_last}
\end{align}
Once the inverse actions of the \(\bS_i\) are available, a new right-hand side is solved by the classical forward sweep
\begin{align}
    \by_1 &= \bS_1^{-1}\bb_1, \\
    \by_i &= \bS_i^{-1}(\bb_i-\bE_i\by_{i-1}),
    \qquad i=2,\ldots,n_b,
    \label{eq:block_forward}
\end{align}
and backward sweep
\begin{align}
    \bu_{n_b} &= \by_{n_b},\\
    \bu_i &= \by_i-\bT_i\bu_{i+1},
    \qquad i=n_b-1,\ldots,1.
    \label{eq:block_backward}
\end{align}
A production classical implementation would normally cache factorizations of the Schur blocks rather than their explicit inverses. Here the explicit inverse columns are organized in a way that they provide small, independent least-squares targets that can be encoded as QUBOs. The decoded approximation \(\widehat{\bS_i^{-1}}\) is assembled once, and {\color{black}the coupling factor \(\widehat{\bT_i}=\widehat{\bS_i^{-1}}\bF_i\) can be cached or applied as the product \(\widehat{\bS_i^{-1}}(\bF_i\bu_{i+1})\) during the backward sweep. The numerical implementation uses the latter equivalent form.} In the five-PDE validation, each subsequent Schur block is formed using the preceding decoded inverse factor. The reported hardware experiment instead uses classically prepared reference Schur blocks.

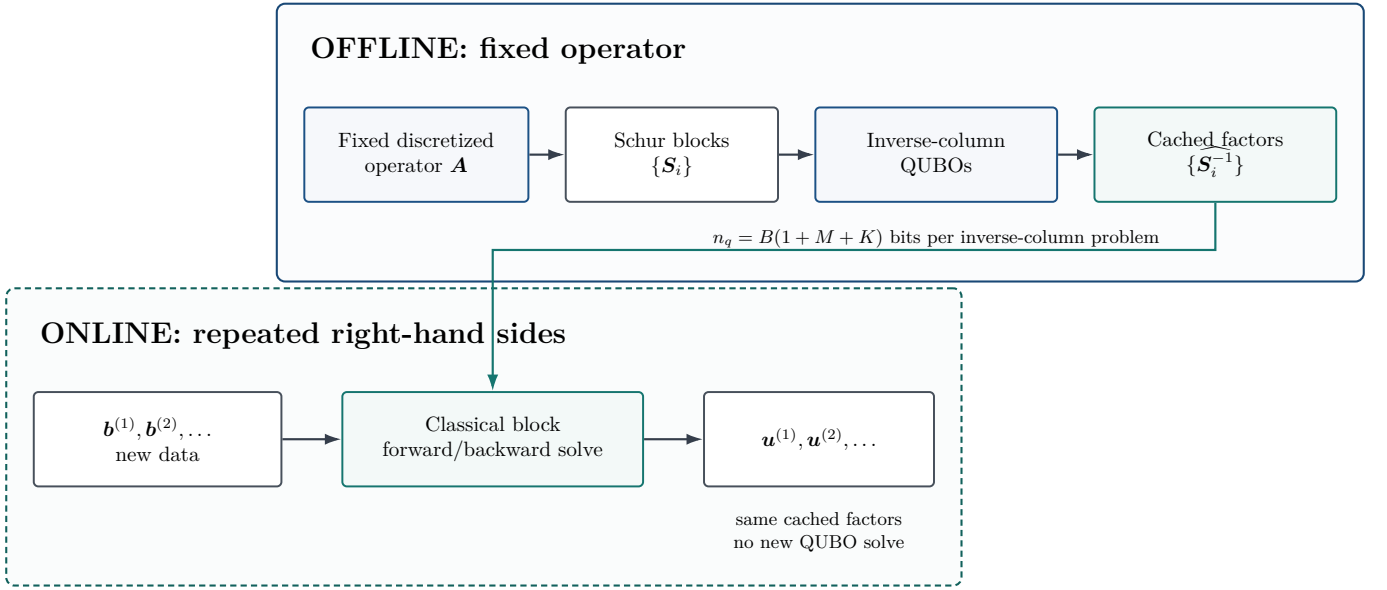
\begin{figure*}[t]
\centering
\resizebox{0.98\textwidth}{!}{%
\begin{tikzpicture}[node distance=8mm and 6mm]
    \tikzset{
      profbox/.style={qpdebox, minimum height=16mm, inner xsep=4mm, inner ysep=3mm, font=\normalsize, line width=0.95pt},
      profaccent/.style={profbox, draw=qpdeblue!90!black, fill=qpdeblue!4},
      profcache/.style={profbox, draw=qpdeteal!90!black, fill=qpdeteal!5},
      profarrow/.style={-{Latex[length=2.8mm,width=1.9mm]}, line width=1.05pt, draw=qpdedark!90}
    }

    \node[anchor=west,font=\Large\bfseries,text=black] (offtitle) at (0,0) {OFFLINE: fixed operator};
    \node[profaccent, below=6mm of offtitle.south west, anchor=north west, text width=30mm,text=black] (A)
      {Fixed discretized\\operator $\bA$};
    \node[profbox, right=of A, text width=28mm,text=black] (schur)
      {Schur blocks\\$\{\bS_i\}$};
    \node[profaccent, right=of schur, text width=33mm,text=black] (qubo)
      {Inverse-column\\QUBOs};
    \node[profcache, right=of qubo, text width=33mm,text=black] (cache)
      {Cached factors\\$\{\widehat{\bS_i^{-1}}\}$};

    \draw[profarrow] (A) -- (schur);
    \draw[profarrow] (schur) -- (qubo);
    \draw[profarrow] (qubo) -- (cache);

    \node[font=\small, text=black, align=center, below=3.2mm of qubo] (bits)
      {$n_q=B(1+M+K)$ bits per inverse-column problem};

    \node[anchor=west,font=\Large\bfseries, below=19mm of A.south west,text=black] (ontitle)
      {ONLINE: repeated right-hand sides};
    \node[profbox, below=6mm of ontitle.south west, anchor=north west, text width=34mm,text=black] (rhs)
      {$\bb^{(1)},\bb^{(2)},\ldots$\\new data};
    \node[profcache, right=10mm of rhs, text width=43mm,text=black] (solve)
      {Classical block\\forward/backward solve};
    \node[profbox, right=10mm of solve, text width=31mm,text=black] (sol)
      {$\bu^{(1)},\bu^{(2)},\ldots$};
    \draw[profarrow] (rhs) -- (solve);
    \draw[profarrow] (solve) -- (sol);
    \draw[profarrow, draw=qpdeteal!90!black] (cache.south) -- ++(0,-8mm) -| (solve.north);
    \node[font=\small, text=black, align=center, below=3mm of sol] (reuse)
      {same cached factors\\no new QUBO solve};

    \begin{scope}[on background layer]
      \node[rounded corners=4pt, draw=qpdeblue!80!black, line width=1.0pt, fill=qpdeblue!2,
            fit=(offtitle)(A)(schur)(qubo)(cache)(bits), inner sep=4.5mm] {};
      \node[rounded corners=4pt, draw=qpdeteal!75!black, densely dashed, line width=1.0pt, fill=qpdeteal!2,
            fit=(ontitle)(rhs)(solve)(sol)(reuse), inner sep=4.5mm] {};
    \end{scope}
\end{tikzpicture}%
}
\caption{Offline--online architecture of the reusable block-Schur QUBO method. For a fixed finite-difference operator, the block-tridiagonal organization is reduced to Schur blocks and each required inverse column is represented as a signed fixed-point least-squares QUBO. The decoded columns are assembled into approximate Schur inverses and cached. \rev{Coupling factors may either be cached explicitly or evaluated from the cached inverse factors during substitution.} Subsequent solves use only the classical block forward/backward recursions in Eqs.~\eqref{eq:block_forward}--\eqref{eq:block_backward}. The largest logical optimization problem therefore contains \(B(1+M+K)\) binary variables and depends on the local block size and numerical precision, not directly on the global number of PDE unknowns \(N\).}
\label{fig:offline_online_architecture}
\end{figure*}

Figure~\ref{fig:offline_online_architecture} shows this separation between operator-dependent preprocessing and repeated classical solves. The QUBO optimizer appears only in the offline stage. Once the cached factors have been constructed, changing the forcing, boundary data, or explicit nonlinear contribution changes only the right-hand side and does not require another QUBO call.

\subsection{Inverse-column QUBO}
\label{sec:inverse_column_qubo}

For a Schur block \(\bS\in\R^{B\times B}\), its \(j\)th inverse column is the solution of
\begin{equation}
    \bS\by_j=\be_j,
\end{equation}
where \(\be_j\) is the \(j\)th Euclidean basis vector. We replace direct inversion by the least-squares objective
\begin{equation}
    E_j(\by)=\|\bS\by-\be_j\|_2^2.
    \label{eq:inverse_column_ls}
\end{equation}
Each component is encoded by one sign bit, \(M\) integer bits, and \(K\) fractional bits,
\begin{equation}
    y_\ell=\gamma\!\left(
      -2^M s_\ell
      +\sum_{m=0}^{M-1}2^m z_{\ell m}
      +\sum_{k=1}^{K}2^{-k}r_{\ell k}
    \right),
    \label{eq:encoding}
\end{equation}
with \(s_\ell,z_{\ell m},r_{\ell k}\in\{0,1\}\) and scale \(\gamma>0\). Hence \(q=1+M+K\) bits encode one scalar and
\begin{equation}
    n_q=Bq=B(1+M+K)
    \label{eq:local_logical_bits}
\end{equation}
encode one inverse column.

Collect all bits into \(\bq\in\{0,1\}^{n_q}\) and write \(\by=\bP\bq\), where \(\bP\) contains the signed fixed-point weights. Substitution in Eq.~\eqref{eq:inverse_column_ls} gives
\begin{align}
    E_j(\bq)
    &=\|\bS\bP\bq-\be_j\|_2^2 \\
    &=\bq^{\mathsf T}\underbrace{\bP^{\mathsf T}\bS^{\mathsf T}\bS\bP}_{\bH}\bq
      -2\underbrace{(\bP^{\mathsf T}\bS^{\mathsf T}\be_j)^{\mathsf T}}_{-\bm{h}^{\mathsf T}/2}\bq+1 .
    \label{eq:qubo_expansion}
\end{align}
Equivalently, with \(\bm{h}=-2\bP^{\mathsf T}\bS^{\mathsf T}\be_j\), an upper-triangular QUBO convention can be written as
\begin{align}
    E_j(\bq) &= 1+\sum_{a\le b}Q_{ab}q_aq_b,\notag\\
    Q_{aa} &= H_{aa}+h_a,\notag\\
    Q_{ab} &= 2H_{ab},\qquad a<b.
    \label{eq:qubo_coefficients}
\end{align}
where \(q_a^2=q_a\) has been used to absorb the linear terms into the diagonal. The additive constant does not affect the minimizer and can be omitted before submission to an optimizer.

{\color{black}For residual correction with target $\bm r$ and scale-dependent encoding $\bP_p=\bP(\gamma_p)$, the corresponding QUBO objective is
\begin{align}
 E_p(\bq)&=\|\bS\bP_p\bq-\bm r\|_2^2 \notag\\
 &=\bq^{\mathsf T}\bP_p^{\mathsf T}\bS^{\mathsf T}\bS\bP_p\bq \notag\\
 &\quad-2\bm r^{\mathsf T}\bS\bP_p\bq+\|\bm r\|_2^2.
 \label{eq:residual_qubo}
\end{align}
The omitted constant is $\|\bm r\|_2^2$, which need not equal one.}

The scale \(\gamma\) controls dynamic range independently of the fractional resolution. It can be selected globally for a fixed operator or separately for each Schur block. In the large-block study of Section~\ref{sec:experiments}, diagonally dominant blocks use the conservative choice
\begin{equation}
    \alpha_\ell^{(i)}
    =|s_{\ell\ell}^{(i)}|-\sum_{m\neq\ell}|s_{\ell m}^{(i)}|,
    \qquad
    \gamma_i=2^{-M}\!\left(\min_\ell\alpha_\ell^{(i)}\right)^{-1},
    \label{eq:gamma_i_scaling}
\end{equation}

when the minimum margin is positive; otherwise, the implementation
uses the inverse-norm fallback
\(\gamma_i=2^{-M}\|\bS_i^{-1}\|_\infty\).
Even when the dominance margin is positive, a very small margin
can produce an excessively large scale and hence coarse absolute
fixed-point resolution. Neither scaling rule alone strictly
guarantees the absence of clipping at the positive representable
endpoint. The representability condition is analyzed in
Section~\ref{sec:analysis}.

\subsection{Offline and online algorithms}

The offline stage consists of assembling \(\bA\), forming the Schur sequence \(\{\bS_i\}\), solving the \(B\) inverse-column QUBOs for every block, and caching \(\widehat{\bS_i^{-1}}\){\color{black}, with optional caching of \(\widehat{\bT_i}\)}. The optimizer used for the column problems is left modular: {\color{black}exact reduced-grid optimization,} classical heuristics, quantum-inspired solvers, or compatible hardware can be substituted without changing the surrounding block algorithm. The online stage then applies Eqs.~\eqref{eq:block_forward}--\eqref{eq:block_backward} for every new \(\bb\). Local inverse residuals \(\|\bS_i\widehat{\bS_i^{-1}}-\bI\|\) provide an operator-level diagnostic before the factors are reused, while the global residual \(\|\bA\widehat{\bu}-\bb\|\) measures the quality of each completed PDE solve.

\section{Error and resource analysis}
\label{sec:analysis}

The binary representation introduces two distinct limitations: finite resolution and finite dynamic range. The former improves systematically with \(K\); the latter must be controlled by \(M\) and \(\gamma\). We first quantify the fixed-point error for one inverse column and then propagate it through the cached block solve. The analysis assumes exact minimization over the chosen finite grid; optimizer error is treated separately in the numerical study.

\subsection{No-clipping condition and fixed-point resolution}

Let \(\by_j^*=\bS^{-1}\be_j\) be an exact inverse column. For the encoding in Eq.~\eqref{eq:encoding}, the scalar grid spacing is
\begin{equation}
    \Delta_K=\gamma 2^{-K},
    \label{eq:grid_spacing}
\end{equation}
and the representable interval is
\begin{equation}
    -2^M\gamma
    \le y_\ell\le
    \gamma(2^M-2^{-K}).
    \label{eq:representable_interval}
\end{equation}
The convergence argument requires that the exact inverse column lie inside this interval, apart from rounding to the nearest grid point. If an entry clips at either boundary, increasing \(K\) alone cannot reduce that saturation error. Thus \(M\) and \(\gamma\) determine dynamic range, whereas \(K\) determines resolution.

Under the no-clipping assumption, there exists a nearest grid vector \(\by_j^\dagger\) satisfying
\begin{equation}
    \|\by_j^\dagger-\by_j^*\|_\infty\le\frac{\Delta_K}{2},
    \qquad
    \|\by_j^\dagger-\by_j^*\|_2
    \le\frac{\sqrt{B}}{2}\Delta_K.
    \label{eq:nearest_grid}
\end{equation}
Let \(\by_{j,K}\) denote an exact minimizer of Eq.~\eqref{eq:inverse_column_ls} over the finite grid. Optimality relative to \(\by_j^\dagger\) yields
\begin{align}
    \|\bS\by_{j,K}-\be_j\|_2
    &\le \|\bS(\by_j^\dagger-\by_j^*)\|_2 \\
    &\le \frac{\|\bS\|_2\sqrt{B}\gamma}{2}\,2^{-K}.
    \label{eq:residual_bound}
\end{align}
Therefore the local inverse-column residual is \(O(2^{-K})\), while its squared least-squares/QUBO energy satisfies
\begin{equation}
    E_{j,K}
    \le
    \frac{\|\bS\|_2^2B\gamma^2}{4}\,4^{-K}.
    \label{eq:energy_bound}
\end{equation}
The distinction between these rates is important experimentally: the decoded residual and inverse error are first order in the grid spacing, whereas the minimized energy is second order.

\subsection{Inverse-factor and PDE-solution error}

Since
\(\bS(\by_{j,K}-\by_j^*)=\bS\by_{j,K}-\be_j\), the smallest singular value of \(\bS\) gives
\begin{align}
    \|\by_{j,K}-\by_j^*\|_2
    &\le \frac{1}{\sigma_{\min}(\bS)}
       \|\bS\by_{j,K}-\be_j\|_2 \\
    &\le \frac{\kappa_2(\bS)\sqrt{B}\gamma}{2}\,2^{-K},
    \label{eq:column_error_bound}
\end{align}
where \(\kappa_2(\bS)=\|\bS\|_2/\sigma_{\min}(\bS)\). Conditioning therefore controls how strongly fixed-point residual error is amplified into inverse-column error.

Let \(\widehat{\bS}^{-1}_K\) be the matrix formed from the \(B\) approximate columns. Summing the columnwise estimate gives
\begin{equation}
    \|\widehat{\bS}^{-1}_K-\bS^{-1}\|_F
    \le
    \frac{\kappa_2(\bS)B\gamma}{2}\,2^{-K},
    \label{eq:inverse_bound}
\end{equation}
so \(\|\widehat{\bS}^{-1}_K-\bS^{-1}\|_2=O(2^{-K})\).

For the complete block solve, denote the exact block-Thomas solution operator by \(\mathcal{B}=\bA^{-1}\) and the cached-QUBO operator by \(\widehat{\mathcal{B}}_K\). Suppose each local inverse satisfies
\begin{equation}
    \|\widehat{\bS}_{i,K}^{-1}-\bS_i^{-1}\|_2
    \le C_i2^{-K}
    \label{eq:local_inverse_assumption}
\end{equation}
and that the exact Schur recursion is stable, with uniformly bounded inverse and coupling factors over the blocks considered.
{\color{black}In the recursively quantized construction, write $\widehat{\bY}_i$ for the decoded approximate inverse and
$\widetilde{\bS}_i=\bD_i-\bE_i\widehat{\bY}_{i-1}\bF_{i-1}$ for $i\ge2$, with $\widetilde{\bS}_1=\bS_1$. Its deviation from the exact Schur block is
\begin{align}
 \bm\Delta_i&=\widetilde{\bS}_i-\bS_i \notag\\
 &=-\bE_i(\widehat{\bY}_{i-1}-\bS_{i-1}^{-1})\bF_{i-1}.
 \label{eq:recursive_schur_perturbation}
\end{align}
If $\theta_i=\|\bS_i^{-1}\|_2\|\bm\Delta_i\|_2<1$, then $\widetilde{\bS}_i$ is invertible and
\begin{equation}
 \|\widetilde{\bS}_i^{-1}-\bS_i^{-1}\|_2
 \le\frac{\|\bS_i^{-1}\|_2^2\|\bm\Delta_i\|_2}{1-\theta_i}.
 \label{eq:recursive_inverse_perturbation}
\end{equation}
Together with the local residual $\bm R_i=\widetilde{\bS}_i\widehat{\bY}_i-\bI$, this yields
\begin{align}
 \|\widehat{\bY}_i-\bS_i^{-1}\|_2
 &\le\|\widetilde{\bS}_i^{-1}\|_2\|\bm R_i\|_2 \notag\\
 &\quad+\|\widetilde{\bS}_i^{-1}-\bS_i^{-1}\|_2.
 \label{eq:recursive_total_inverse_error}
\end{align}
Thus the stated local inverse assumption requires control of both the quantized inverse residual and the propagated Schur perturbation; it does not follow from local grid resolution alone.} Standard perturbation of the finite composition of forward and backward block maps then gives
\begin{equation}
    \|\widehat{\mathcal{B}}_K-\mathcal{B}\|_2
    \le C_{\mathrm{BT}}2^{-K}+O(2^{-2K}),
    \label{eq:block_operator_error}
\end{equation}
where \(C_{\mathrm{BT}}\) depends on the block couplings, the number of blocks \(n_b\), and the norms and condition numbers of the Schur complements. In particular, perturbations accumulate through the forward and backward recursions, so \(C_{\mathrm{BT}}\) need not remain bounded as \(n_b\) increases; its growth depends on the products of the corresponding coupling operators. Hence, for a fixed right-hand side,
\begin{equation}
    \|\widehat{\bu}_K-\bu\|_2
    \le \left(C_{\mathrm{BT}}2^{-K}+O(2^{-2K})\right)\|\bb\|_2.
    \label{eq:pde_error_bound}
\end{equation}

For repeated time stepping, write the right-hand side as
\(\bb^n=\bm{G}(\bu^n)\), and suppose \(\bm{G}\) is Lipschitz
with constant \(L_G\). If \(e_n=\|\widehat{\bu}^{\,n}-\bu^n\|_2\), then
\begin{equation}
    e_{n+1}
    \le
    \|\widehat{\mathcal{B}}_K\|_2 L_G e_n
    +
    \|\widehat{\mathcal{B}}_K-\mathcal{B}\|_2
    \|\bm{G}(\bu^n)\|_2.
\end{equation}
Thus, local inverse-factor errors may accumulate or be amplified
over repeated time steps, depending on the stability of the
time-stepping map.

Equations~\eqref{eq:energy_bound}, \eqref{eq:inverse_bound}, and \eqref{eq:pde_error_bound} therefore predict \(2^{-K}\) convergence for inverse/solution errors and \(4^{-K}\) convergence for the minimized QUBO energy, provided that clipping and optimizer error are absent.

\subsection{Logical variables, couplings, and amortization}
\label{sec:resource_analysis}

Let \(q=1+M+K\) denote the number of bits per real scalar. A direct fixed-point encoding of the full \(N\)-unknown system would expose
\begin{equation}
    n_{\mathrm{full}}=Nq
\end{equation}
logical binary variables to one global QUBO. In contrast, Eq.~\eqref{eq:local_logical_bits} limits every inverse-column QUBO to
\begin{equation}
    n_{\mathrm{local}}=Bq,
    \qquad
    \frac{n_{\mathrm{local}}}{n_{\mathrm{full}}}=\frac{B}{N}.
    \label{eq:logical_reduction}
\end{equation}
This reduction concerns the largest logical QUBO dimension, not the total amount of offline work. There are \(n_bB=N\) inverse columns across all Schur blocks, so one complete cached factorization requires \(N\) local QUBO solves before any optional residual-refinement passes.

The least-squares construction can be logically dense even when the underlying PDE matrix is sparse. With an upper-triangular storage convention, a fully dense local QUBO contains
\begin{equation}
    m_Q^{\max}
    =\frac{n_{\mathrm{local}}(n_{\mathrm{local}}+1)}{2}
    =O(B^2q^2)
    \label{eq:coupling_scaling}
\end{equation}
coefficients. We therefore report the logical density
\begin{equation}
    \rho_Q=
    \frac{\#\{(a,b):Q_{ab}\neq0,\ a\le b\}}
         {n_{\mathrm{local}}(n_{\mathrm{local}}+1)/2}.
    \label{eq:qubo_density}
\end{equation}
The distinction between variable count and coupling count is essential for hardware interpretation because embedding, routing, coefficient precision, or native polynomial connectivity can depend strongly on the interaction graph.

If exhaustive enumeration is used as an exact finite-bit oracle, the offline search work scales as
\begin{equation}
    O\!\left(N2^{Bq}\right),
    \label{eq:exact_oracle_complexity}
\end{equation}
which is exponential in the local encoded dimension. {\color{black}The present $B=2$ experiments instead use an exact reduced-grid minimizer (validated against exhaustive search for small instances); neither method is the proposed scalable optimizer.} For \(R\) subsequent right-hand sides, the factor-construction cost is paid once and its average contribution per solve is reduced by \(1/R\); the online stage remains the classical block substitution of Eqs.~\eqref{eq:block_forward}--\eqref{eq:block_backward}. No quantum or heuristic speedup follows from Eq.~\eqref{eq:logical_reduction} alone.

\section{Numerical and hardware protocol}
\label{sec:experiments}

The experiments separate quantization error, representability, heuristic optimizer error, and hardware sampling into distinct tests. Table~\ref{tab:benchmark_protocol} summarizes the five finite-difference benchmark families. The exact-oracle studies establish the finite-bit theory at small local dimension; the large-block study tests whether the same representation remains usable as \(B\) grows; the \(B=8\) study probes non-exact optimization at 96 bits; and the hardware experiment executes 24-bit PDE-derived QUBOs on QCI Dirac-3.

\begin{table*}[t]
\centering
\small
\caption{Finite-difference benchmark families and their role in the validation protocol. The fixed-validation systems all use \(B=2\), \(M=1\), and \(K=10\); the parameters listed below specify the physical/discretization setting. The Poisson family is additionally used for the \(N=960\) large-block experiments, and Heat, Poisson, and Klein--Gordon are reused in the Dirac-3 hardware campaign. Manufactured/discrete reference solutions are used where indicated so that solver error can be separated from finite-difference truncation error.}
\label{tab:benchmark_protocol}
\begin{tabularx}{\textwidth}{@{}p{2.25cm}p{2.9cm}p{2.2cm}X@{}}
\toprule
Benchmark & Fixed operator / character & Fixed-validation size & Boundary conditions and principal parameters \\
\midrule
Heat 1D & Implicit diffusion; parabolic & \(N=20\) & Time-dependent Dirichlet data; {\color{black}\(\nu=0.08\),} \(\Delta t=10^{-3}\), \(T=0.05\). \\
Burgers 1D & Implicit diffusion + explicit convection; nonlinear IMEX & \(N=20\) & Dirichlet data; \(\nu=0.04\), \(\Delta t=10^{-3}\), \(T=0.05\). \\
Klein--Gordon 1D & Coupled first-order wave/reaction system; hyperbolic & state dimension \(N=40\) & Interleaved \([u_1,v_1,u_2,v_2,\ldots]^\mathsf T\) ordering; Dirichlet data; {\color{black}\(c=1.1\), \(\mu=2\),} \(\Delta t=10^{-3}\), \(T=0.05\). \\
Poisson 2D & Steady elliptic & \(10\times2\) interior grid & Dirichlet manufactured data; {\color{black}\(\Delta x=1/11\), \(\Delta y=1/3\); }\(\bb=\bA\bu_{\mathrm{exact}}\) in the fixed validation to isolate cached-factor error. \\
Helmholtz 2D & Steady shifted elliptic/oscillatory & \(10\times2\) interior grid & Dirichlet manufactured data; algebraic operator \(-\Delta-k^2\bI\), \(k=3\), {\color{black}\(\Delta x=1/11\), \(\Delta y=1/3\).} \\
\bottomrule
\end{tabularx}
\end{table*}

{\color{black}The five-PDE validation imposes the Dirichlet traces of the following manufactured fields (with $x,y\in[0,1]$):
\begin{align*}
 u_{\mathrm H}^*(x,t)&=e^{-t}[\sin(\pi x)+x],\\
 u_{\mathrm B}^*(x,t)&=0.35+0.15x+0.25e^{-t}\sin(\pi x),\\
 u_{\mathrm{KG}}^*(x,t)&=\cos(1.7t)[\sin(\pi x)+0.2x],\\
 u_{\mathrm P}^*(x,y)&=\sin(\pi x)\sin(2\pi y)
 +0.15x+0.3y+0.15xy,\\
 u_{\mathrm{Helm}}^*(x,y)&=\cos(2\pi x)\sin(\pi y)
 +0.1x+0.05y.
\end{align*}
For Klein--Gordon, the manufactured velocity is $v_{\mathrm{KG}}^*=\partial_t u_{\mathrm{KG}}^*$. Heat and Burgers include their manufactured forcing and boundary terms at the time levels of the stated implicit and IMEX schemes. For the two elliptic systems, the tested discrete right-hand side is set to $\bb=\bA\bu_h^*$; the boundary trace is included in the manufactured field. This isolates the cached-solver error, rather than asserting equality with a right-hand side obtained by sampling the continuous forcing and then applying a separate boundary contribution.}

\subsection{Exact finite-bit oracle and five-PDE validation}

The convergence experiment uses \(B=2\), \(M=1\), and \(K=4,5,\ldots,10\). Each inverse-column problem therefore contains \(n_q=2(K+2)\) bits, ranging from 12 to 24, and is minimized exactly over its finite grid. The two sweep systems are Poisson 2D and Klein--Gordon 1D, with fixed scales \(\gamma=0.01\) and \(\gamma=1.0\), respectively. For the \(N=40\) Klein--Gordon validation, the two state components are interleaved as \([u_1,v_1,u_2,v_2,\ldots]^{\mathsf T}\). The exact finite-grid minimizers isolate the convergence predicted by Section~\ref{sec:analysis}.

{\color{black}The Figure~\ref{fig:k_sweep} precision sweep uses smaller systems than the five-PDE validation of Table~\ref{tab:benchmark_protocol}. Its Poisson case has a $6\times2$ interior grid on $(0,1)^2$, with $h_x=1/7$, $h_y=1/3$, homogeneous Dirichlet data, and discrete reference field
\(u_h^*(x,y)=[\sin(\pi x)+0.2\sin(2\pi x)]\sin(\pi y)\); the right-hand side is set to $\bb=\bA\bu_h^*$. Its Klein--Gordon case uses six interior spatial points ($N=12$ interleaved state variables), $c=1$, $\mu=1.25$, $\Delta t=0.015$, and $T=0.3$ (20 implicit steps), with homogeneous Dirichlet displacement, zero forcing, and initial data $u(x,0)=\sin(\pi x)$, $v(x,0)=0$. For each $K$, both cases use one exact finite-grid optimization per inverse column and recursively form subsequent Schur blocks from the quantized inverse factors, without residual-refinement passes. The Poisson solution is compared with its exact discrete manufactured reference; the Klein--Gordon solution error is measured against the corresponding dense numerical trajectory.}

All five PDE families in Table~\ref{tab:benchmark_protocol} are then tested with
\begin{equation}
    B=2,\qquad M=1,\qquad K=10,
    \label{eq:fixed_validation_encoding}
\end{equation}
corresponding to 24 logical bits per inverse-column QUBO. The initial scales are
$\gamma_0=0.60$ for Heat and Burgers, $\gamma_0=1.00$ for Klein--Gordon,
and $\gamma_0=0.01$ for Poisson and Helmholtz. Each inverse column uses one
initial QUBO solve and three QUBO residual corrections{\color{black}, each using Eq.~\eqref{eq:residual_qubo} with its stage-specific residual and scale}.
For Schur block $\bS_i$ and basis vector $\be_j$, initialize
$\by_{ij}^{(0)}=\bm{0}$ and prescribe
\begin{align}
    \gamma_p &= \gamma_0 8^{-p}, \quad p=0,1,2,3,
        \label{eq:five_pde_scale_schedule}\\
    \bm r_{ij}^{(p)} &= \be_j-\bS_i\by_{ij}^{(p)}, \notag\\
    \bq_{ij,p}^{*} &\in \underset{\bq\in\{0,1\}^{24}}{\operatorname{arg\,min}}\,
      \bigl\|\bS_i\bP(\gamma_p)\bq-\bm r_{ij}^{(p)}\bigr\|_2^2,
      \notag\\
    \by_{ij}^{(p+1)} &= \by_{ij}^{(p)}+\bP(\gamma_p)\bq_{ij,p}^{*},
        \label{eq:five_pde_multiscale}\\
    \textcolor{black}{\bm r_{ij}^{(p+1)}} &= \textcolor{black}{\bm r_{ij}^{(p)}-\bS_i\bP(\gamma_p)\bq_{ij,p}^{*}.}
        \label{eq:five_pde_residual_update}
\end{align}
The shrink factor eight is specified before optimization, not selected from a
continuous correction or the desired numerical result. The multiscale experiment and
its fixed-grid control use the same PDE operators, source terms, time steps,
initial data, manufactured fields and Dirichlet traces; subsequent Schur blocks
use previously decoded quantized inverse factors. A continuous correction is
evaluated only afterwards to audit representability, not to select the scale
or a binary solution.

All finite-bit QUBO minima are obtained with an exact reduced-grid oracle for
two encoded scalars. The oracle minimizes the 24-bit QUBO objective without
visiting all $2^{24}$ binary assignments. The fixed-grid control uses
$\gamma_p=\gamma_0$ at all four stages; its three correction passes give zero updates.

\subsection{Large-block representability and 96-bit heuristic optimization}

To separate block-size scaling from global grid size, a Poisson line-block problem is constructed at fixed \(N=960\) with
\begin{equation}
    B\in\{4,8,16,32\},\qquad M=1,\qquad K=10.
    \label{eq:largeB_protocol}
\end{equation}
The interior grid is \(B\times(960/B)\) and each line of \(B\) unknowns forms one block. {\color{black}It uses the unscaled, equal-spacing five-point stencil with zero exterior (homogeneous Dirichlet) values and the prescribed load $b_{ij}=\sin(\pi i/(n_b+1))\sin(\pi j/(B+1))$, where $n_b=960/B$. As $B$ varies, these rectangular index grids are not interpreted as different discretizations of the same unit-square continuum domain.} The five-point operator is sparse, but the recursively formed Schur complements are dense. For this experiment the exact dense Schur inverse columns are rounded onto the same signed fixed-point grid, using the block-dependent scale in Eq.~\eqref{eq:gamma_i_scaling}, and then propagated through the cached block solve. The experiment therefore tests dynamic range, saturation, conditioning, and resource growth.

A separate \(B=8\) study examines practical optimization of the resulting 96-bit, fully dense inverse-column QUBOs. Thirty-two representative instances are used: all eight columns from Schur blocks 1, 30, 60, and 120. Four strategies are compared on the identical QUBO objective: direct fixed-point rounding, greedy single-bit descent initialized by the rounded state, cold random-start simulated annealing followed by greedy polishing, and warm-start simulated annealing followed by greedy polishing. In addition, greedy descent from the rounded state is applied to all 960 inverse columns to test the downstream effect of QUBO-space polishing on a complete cached factorization. No exact 96-bit ground-state claim is made.

\subsection{Dirac-3 hardware protocol}
\label{sec:dirac_protocol}
QCI Dirac-3 is a qudit entropy quantum device implemented as a hybrid analog optimization system combining quantum optics and digital electronics~\cite{qci_dirac3_beginner_guide,qci_dirac3_spec_sheet,nguyen2024open_quantum_system}. In this experiment, its discrete-Hamiltonian interface is used with binary-valued variables to solve the inverse-column QUBOs. The hardware campaign uses Heat 1D, Poisson 2D, and Klein--Gordon 1D systems whose local inverse-column problems retain \(B=2\), \(M=1\), and \(K=10\). The Heat and Poisson systems contain 20 algebraic unknowns; the hardware Klein--Gordon case uses 10 spatial points in the interleaved two-component state, also giving 20 unknowns. Consequently, each mapped inverse-column instance is a 24-bit QUBO. The upper-triangular least-squares construction contains 300 coefficients and is fully dense (\(\rho_Q=1\)). Before submission, the constant term is removed and the remaining coefficients are divided by the largest absolute nonzero coefficient. 

There are 60 mapped PDE inverse-column instances across the three systems. Repeated Schur/QUBO structures reduce these to 34 unique hardware submissions, whose returned samples are mapped back to the 60 instances before the cached PDE factors are assembled. Thirty-two submissions use 20 returned samples with relaxation schedule 1; two diagnostic submissions use 100 samples with relaxation schedule 4. For each unique job, the lowest recomputed submitted objective among the returned samples is retained. 

{\color{black}The terminal-system benchmarks represented in Figure~\ref{fig:dirac_reconstruction} admit the following homogeneous-Dirichlet finite-difference specifications. The Heat problem is a single backward-Euler step of $u_t=0.08u_{xx}$ on $x_i=i/21$ ($i=1,\ldots,20$), with $u(0,t)=u(1,t)=0$, $u_i^0=\sin(\pi i/21)$ and $\Delta t=10^{-3}$. Its algebraic right-hand side is $b_i=\sin(\pi i/21)$. The Poisson problem is $-\Delta_h u=f_h$ on the $10\times2$ interior grid of $(0,1)^2$, with $h_x=1/11$, $h_y=1/3$, $u|_{\partial\Omega}=0$, and $b_{ij}=f_{h,ij}=\sin(\pi i/11)\sin(\pi j/3)$. The Klein--Gordon problem is the first-order system $u_t=v$, $v_t=1.1^2u_{xx}-2^2u$, discretized by one implicit step with $\Delta t=10^{-3}$ on $x_i=i/11$ ($i=1,\ldots,10$), $u(0,t)=u(1,t)=0$, $u_i^0=\sin(\pi i/11)$ and $v_i^0=0$. Its interleaved algebraic right-hand side is $\bb_i=[\sin(\pi i/11),0]^{\mathsf T}$.
 Thus Figure~\ref{fig:dirac_reconstruction} tests a terminal algebraic solve with one-stage hardware-derived inverse columns.}

The selected Dirac-3 bitstring is compared against the exact 24-bit finite-grid reference for the same normalized submitted objective and is independently decoded back to an inverse column before the PDE reconstruction. The principal error measures used throughout the numerical sections are
\begin{align}
    e_u &= \frac{\|\widehat{\bu}-\bu_{\mathrm{ref}}\|_2}{\|\bu_{\mathrm{ref}}\|_2},\notag\\
    r_A &= \frac{\|\bA\widehat{\bu}-\bb\|_2}{\|\bb\|_2},\notag\\
    r_S &= \|\bS_i\widehat{\by}_{ij}-\be_j\|_2.
    \label{eq:primary_metrics}
\end{align}
For hardware instance \(j\), we additionally report the submitted-objective gap relative to the exact finite-bit optimum,
\begin{equation}
    \Delta E_j=E_{\mathrm{D}}^{(j)}-E_*^{(j)},
    \qquad
    g_j=\frac{\Delta E_j}{\max(|E_*^{(j)}|,10^{-12})},
    \label{eq:hardware_gap}
\end{equation}
so that objective quality, decoded inverse quality, and final PDE accuracy can be assessed separately.

\section{Results}
\label{sec:results}

The experiments are reported in the same order as the validation hierarchy in Section~\ref{sec:experiments}. The first two studies use exact finite-bit minimizers and therefore isolate the consequences of the QUBO representation itself. The large-block study then separates representability from optimization, the 96-bit study probes practical heuristic search, and the final study replaces the exact oracle by samples returned from QCI Dirac-3 hardware.

\subsection{Quantization-controlled convergence}
\label{sec:results_ksweep}

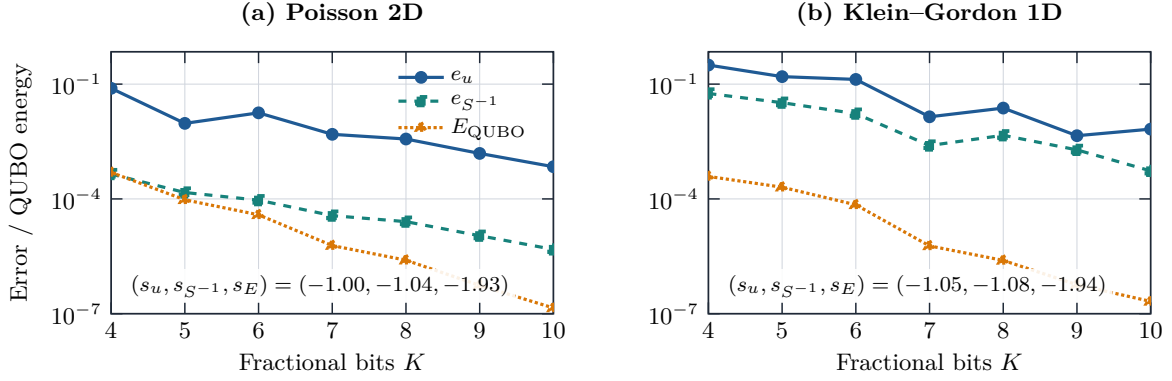
\begin{figure*}[t]
\centering
\begin{tikzpicture}
\begin{groupplot}[
    group style={group size=2 by 1, horizontal sep=2.05cm},
    qpde axis,
    width=0.405\textwidth,
    height=0.275\textwidth,
    ymode=log,
    ymin=1e-7,ymax=7e-1,
    xlabel={Fractional bits $K$},
    xmin=4,xmax=10,
    xtick={4,5,6,7,8,9,10},
    tick label style={font=\small},
    label style={font=\small},
    title style={font=\small\bfseries},
    legend style={font=\footnotesize,draw=none,fill=none},
]
\nextgroupplot[
    title={(a) Poisson 2D},
    ylabel={\rev{Error / QUBO energy}},
    legend pos=north east,
]
\addplot[qpde blue,solid,mark=*] table[col sep=comma,x=K,y=poisson_solution]{k_sweep_compact.csv};
\addlegendentry{\rev{$e_u$}}
\addplot[qpde teal,dashed,mark=square*] table[col sep=comma,x=K,y=poisson_inverse]{k_sweep_compact.csv};
\addlegendentry{\rev{$e_{S^{-1}}$}}
\addplot[qpde orange,densely dotted,mark=triangle*] table[col sep=comma,x=K,y=poisson_energy]{k_sweep_compact.csv};
\addlegendentry{\rev{$E_{\rm QUBO}$}}
\node[anchor=south west,fill=white,fill opacity=0.90,text opacity=1,rounded corners=1pt,inner sep=2.2pt,font=\footnotesize] at (rel axis cs:0.03,0.04)
  {\rev{$(s_u,s_{S^{-1}},s_E)=(-1.00,-1.04,-1.93)$}};

\nextgroupplot[
    title={(b) Klein--Gordon 1D},
]
\addplot[qpde blue,solid,mark=*] table[col sep=comma,x=K,y=kg_solution]{k_sweep_compact.csv};
\addplot[qpde teal,dashed,mark=square*] table[col sep=comma,x=K,y=kg_inverse]{k_sweep_compact.csv};
\addplot[qpde orange,densely dotted,mark=triangle*] table[col sep=comma,x=K,y=kg_energy]{k_sweep_compact.csv};
\node[anchor=south west,fill=white,fill opacity=0.90,text opacity=1,rounded corners=1pt,inner sep=2.2pt,font=\footnotesize] at (rel axis cs:0.03,0.04)
  {\rev{$(s_u,s_{S^{-1}},s_E)=(-1.05,-1.08,-1.94)$}};
\end{groupplot}
\end{tikzpicture}
\caption{Exact-oracle validation of the fixed-point convergence theory for $B=2$, $M=1$, and $K=4,\ldots,10$. Each plotted point uses the exact minimizer of every finite-bit inverse-column QUBO, {\color{black}obtained by exact reduced-grid minimization, so the curves isolate quantization and propagation error rather than heuristic optimizer error.} \rev{Here $e_u$ is the relative PDE solution error, $e_{S^{-1}}$ is the mean Schur-inverse error, and $E_{\rm QUBO}$ is the mean minimized least-squares/QUBO energy.} In both PDEs, the relative solution error and mean Schur-inverse error have fitted $\log_2$ slopes close to $-1$, consistent with $O(2^{-K})$, while the mean minimized least-squares/QUBO energy has a slope close to $-2$, consistent with $O(4^{-K})$. The occasional non-monotonic final-state point reflects propagation of many local quantization perturbations and does not alter the fitted asymptotic trend.}
\label{fig:k_sweep}
\end{figure*}

Figure~\ref{fig:k_sweep} tests the rates derived in Section~\ref{sec:analysis}. For each value of $K$, every inverse-column QUBO is solved by {\color{black}exact reduced-grid minimization, so heuristic optimizer failure is removed from the experiment.} A least-squares fit of $\log_2$ error against $K$, using all seven values $K=4,\ldots,10$, gives solution-error, mean inverse-error, and mean minimized-energy slopes of $-0.997$, $-1.035$, and $-1.932$ for Poisson 2D, and $-1.051$, $-1.080$, and $-1.940$ for Klein--Gordon 1D. These values are close to the predicted slopes $-1$, $-1$, and $-2$, respectively.

For Poisson, the relative solution error decreases from $7.77\times10^{-2}$ at $K=4$ to $6.94\times10^{-4}$ at $K=10$. For Klein--Gordon, the relative final-state error against the dense reference decreases from $3.15\times10^{-1}$ to $6.67\times10^{-3}$ over the same range, with a smaller value $4.49\times10^{-3}$ at $K=9$. The modest non-monotonicity of the final PDE error at individual $K$ values does not contradict the local estimate: quantization errors are propagated through many block substitutions and, for Klein--Gordon, many time steps, so different local perturbations can combine constructively or destructively. The fitted trend nevertheless follows the predicted first-order decay, while the squared QUBO residual follows the predicted second-order decay.

% layout barrier relaxed

\subsection{Five-PDE validation with a common 24-bit local encoding}
\label{sec:results_fixed}

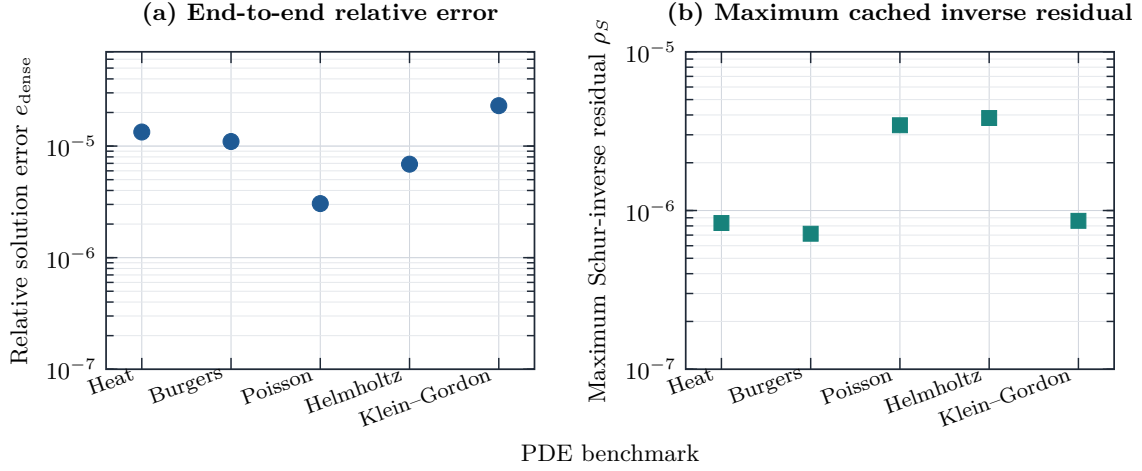
\begin{figure*}[t]
\centering
\begin{tikzpicture}
\begin{groupplot}[
    group style={group size=2 by 1, horizontal sep=2.00cm},
    qpde axis,
    width=0.395\textwidth,
    height=0.315\textwidth,
    ymode=log,
    x tick label style={rotate=20,anchor=east,font=\footnotesize},
    symbolic x coords={Heat,Burgers,Poisson,Helmholtz,Klein--Gordon},
    xtick=data,
    enlarge x limits=0.10,
    ymin=1e-7,
    tick label style={font=\small},
]
\nextgroupplot[title={(a) End-to-end relative error},ylabel={\rev{Relative solution error $e_{\rm dense}$}},ymax=7e-5]
\addplot[only marks,mark=*,mark size=3.0pt,color=qpdeblue]
  table[col sep=comma,x=short,y=rel_dense]{figure3_multiscale_plot.csv};
\nextgroupplot[title={(b) Maximum cached inverse residual},ylabel={\rev{Maximum Schur-inverse residual $\rho_S$}},ymax=1e-5]
\addplot[only marks,mark=square*,mark size=2.8pt,color=qpdeteal]
  table[col sep=comma,x=short,y=max_inv_res]{figure3_multiscale_plot.csv};
\end{groupplot}
\node[font=\small] at ($(group c1r1.south)!0.5!(group c2r1.south)+(0,-11mm)$) {\rev{PDE benchmark}};
\end{tikzpicture}
\caption{Exact finite-bit five-PDE validation using the prescribed
multiscale correction schedule $\gamma_p=\gamma_0/8^p$, $p=0,1,2,3$.
Panel (a) reports relative error against the dense numerical solve, and
panel (b) reports the maximum cached Schur-inverse residual. Each point
represents the complete cached-factor PDE calculation with the 24-bit
inverse-column encoding.}
\label{fig:fixed_validation}
\end{figure*}

The fixed-parameter study tests whether the same local encoding can be reused
across qualitatively different PDE operators, with the prescribed multiscale
corrections of Eq.~\eqref{eq:five_pde_scale_schedule}. Table~\ref{tab:fixed_validation}
reports end-to-end accuracy, conditioning, local inverse residuals and QUBO
workload. Figure~\ref{fig:fixed_validation} presents the two principal
accuracy diagnostics.

\begin{table*}[t]
\centering
\scriptsize
\caption{Exact finite-bit multiscale validation for all five PDE families,
with $B=2$, $M=1$, $K=10$ (24 logical variables per inverse-column QUBO),
one initial solve and three residual corrections using
$\gamma_p=\gamma_0/8^p$ for $p=0,1,2,3$. $e_{\rm dense}$ is relative
$\ell_2$ error against the dense numerical reference; $r_\infty$ is the final
absolute infinity-norm algebraic residual, evaluated with the actual last-step
right-hand side for time-dependent problems; and
$\rho_S=\max_i\|{\color{black}\widetilde{\bS}_i\widehat{\bY}_i}-\bI\|_F$. {\color{black}Here $\widetilde{\bS}_i$ is the recursively constructed Schur block and $\widehat{\bY}_i$ is its decoded inverse approximation.} The condition estimates
refer to the assembled matrix and the most ill-conditioned {\color{black}recursively constructed} local Schur block.
``Calls / updates'' counts exact QUBO solves and nonzero residual corrections.
All QUBO minima are obtained with the exact reduced-grid oracle.}
\label{tab:fixed_validation}
\begin{tabular}{lrrrrrrrr}
\toprule
PDE & $\gamma_0$ & $N$ & $\kappa_2(A)$ & $e_{\rm dense}$ & $r_\infty$ & $\rho_S$ & $\max\kappa_2({\color{black}\widetilde{\bS}_i})$ & Calls / updates \\
\midrule
Heat 1D & 0.60 & 20 & 1.14 & $1.34\times10^{-5}$ & $7.12\times10^{-7}$ & $8.34\times10^{-7}$ & 1.07 & 80 / 60 \\
Burgers 1D & 0.60 & 20 & 1.07 & $1.10\times10^{-5}$ & $2.73\times10^{-7}$ & $7.14\times10^{-7}$ & 1.03 & 80 / 60 \\
Poisson 2D & 0.01 & 20 & 26.66 & $3.05\times10^{-6}$ & $3.04\times10^{-4}$ & $3.45\times10^{-6}$ & 1.21 & 80 / 60 \\
Helmholtz 2D & 0.01 & 20 & 50.21 & $6.88\times10^{-6}$ & $1.67\times10^{-4}$ & $3.83\times10^{-6}$ & 1.33 & 80 / 60 \\
Klein--Gordon 1D & 1.00 & 40 & 6.35 & $2.30\times10^{-5}$ & $2.75\times10^{-7}$ & $8.61\times10^{-7}$ & 2.79 & 160 / 120 \\
\bottomrule
\end{tabular}
\end{table*}

Heat and Burgers yield relative errors of $1.34\times10^{-5}$ and
$1.10\times10^{-5}$ against separate dense numerical solutions. The Burgers
case retains nonlinear IMEX evolution with a fixed implicit diffusion
operator; only the right-hand side changes, so the same cached inverse
factors are reused at every step. The prescribed scale schedule improves
resolution of residual corrections within the 24-bit QUBO formulation.

The Poisson and Helmholtz assembled operators have condition estimates $26.66$
and $50.21$, while their maximum local Schur condition numbers remain $1.21$
and $1.33$. Their relative errors are $3.05\times10^{-6}$ and
$6.88\times10^{-6}$. Both use the manufactured discrete field as the dense
reference; these errors quantify cached-factor approximation rather than
finite-difference truncation error. The larger absolute elliptic
$r_\infty$ values are not directly comparable to the relative solution errors.

For Klein--Gordon, the interleaved state has $N=40$ and requires 160 QUBO
calls. Its combined-state relative error is $2.30\times10^{-5}$; the
separate displacement and velocity relative errors are
$1.42\times10^{-5}$ and $1.28\times10^{-4}$, respectively. The terminal
state correlation with the dense numerical reference is approximately
$0.99999999985$.

The multiscale experiment uses 80 QUBO calls per $N=20$ benchmark and 160
for Klein--Gordon, for 480 exact finite-grid minimizations in total. All 360
correction solves give nonzero updates; all retrospectively checked continuous
corrections are representable on their prescribed stage grids. The fixed-grid
control, in contrast, returns zero for all 360 correction solves. Each QUBO has $2^{24}$ logical assignments;
the reduced-grid solver finds its minimum without visiting the full binary
state space.

% layout barrier relaxed

\subsection{Large-block representability and resource scaling}
\label{sec:results_largeB}

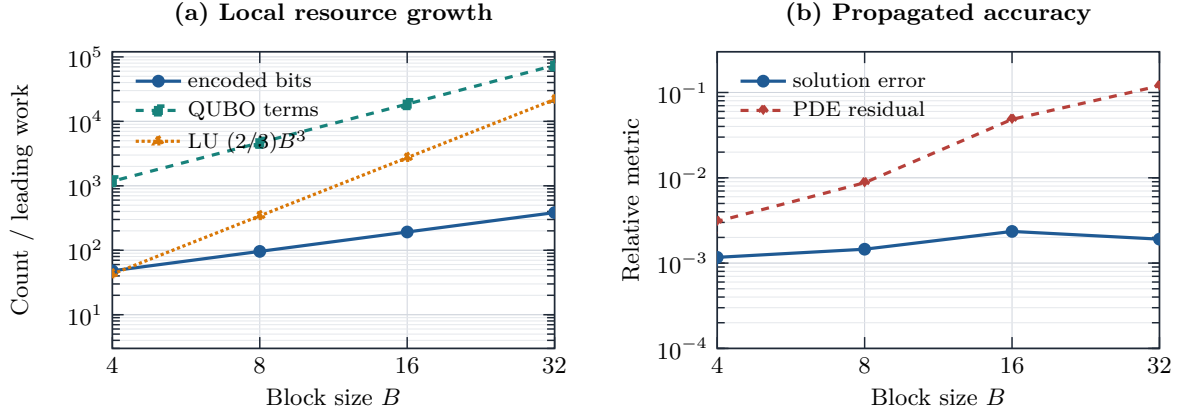
\begin{figure*}[t]
\centering
\begin{tikzpicture}
\begin{groupplot}[
    group style={group size=2 by 1, horizontal sep=2.15cm},
    qpde axis,
    width=0.405\textwidth,
    height=0.30\textwidth,
    xmode=log,log basis x=2,
    xmin=4,
    xtick={4,8,16,32},
    xticklabels={4,8,16,32},
    xlabel={Block size $B$},
]
\nextgroupplot[title={(a) Local resource growth},ymode=log,ylabel={Count / leading work},ymin=3,ymax=1.2e5,legend pos=north west]
\addplot[qpde blue,solid,mark=*] table[col sep=comma,x=B,y=encoded_bits_per_inverse_column]{largeB.csv};
\addlegendentry{encoded bits}
\addplot[qpde teal,dashed,mark=square*] table[col sep=comma,x=B,y=upper_triangular_qubo_terms]{largeB.csv};
\addlegendentry{QUBO terms}
\addplot[qpde orange,densely dotted,mark=triangle*] table[col sep=comma,x=B,y=dense_inverse_flops_2over3B3]{largeB.csv};
\addlegendentry{{\color{black}LU }$(2/3)B^3$}
\nextgroupplot[title={(b) Propagated accuracy},ymode=log,ylabel={Relative metric},ymin=1e-4,ymax=3e-1,legend pos=north west]
\addplot[qpde blue,solid,mark=*] table[col sep=comma,x=B,y=relative_error_vs_dense]{largeB.csv};
\addlegendentry{solution error}
\addplot[qpde red,dashed,mark=diamond*] table[col sep=comma,x=B,y=relative_residual]{largeB.csv};
\addlegendentry{PDE residual}
\end{groupplot}
\end{tikzpicture}
\caption{Large-block Poisson line-block study at fixed $N=960$. Panel (a) separates three different scaling quantities: the inverse-column logical variable count grows linearly as $12B$, the dense upper-triangular QUBO coefficient count grows quadratically, and {\color{black}the leading LU-factorization work for a conventional dense $B\times B$ block is $(2/3)B^3$ and grows cubically; explicitly forming the inverse incurs additional work.} Panel (b) shows the solution error and global residual after rounded fixed-point Schur inverses are cached and used in the complete block solve. Per-block scaling produces zero clipping for all four block sizes (Table~\ref{tab:largeB}); the solution remains accurate at the $10^{-3}$ level while the residual becomes more sensitive as $B$ increases.}
\label{fig:largeB}
\end{figure*}

The $N=960$ Poisson line-block experiment addresses a different question from the exact-oracle tests: whether the local fixed-point representation remains numerically usable when the block contains many unknowns. Table~\ref{tab:largeB} and Figure~\ref{fig:largeB} summarize the sweep from $B=4$ to $B=32$. The encoded inverse-column dimension grows exactly as $12B$, from 48 to 384 bits. Over the same interval the number of dense upper-triangular QUBO coefficients grows from 1176 to 73,920, while the {\color{black}leading dense LU-factorization work} $(2/3)B^3$ grows from 42.67 to $2.18\times10^4$ floating-point operations per block {\color{black}factorization}. Thus the logical variable count is linear in $B$ at fixed precision, but the interaction count is quadratic and the corresponding classical {\color{black}dense LU-factorization work} is cubic.

\begin{table*}[t]
\centering
\scriptsize
\caption{Large-$B$ Poisson line-block representability study at fixed global dimension $N=960$, $M=1$, and $K=10$. Each line of $B$ grid unknowns forms one block. The exact dense Schur inverse columns are rounded onto the same fixed-point grid that would be used by the QUBO and then propagated through the cached block solve; therefore the table measures representation, scaling, and propagation, not successful global optimization of the 96--384-bit QUBOs. Saturation is the fraction of rounded inverse entries clipped at the representable bounds. {\color{black}The $(2/3)B^3$ column gives the leading dense LU-factorization work per block.}}

\label{tab:largeB}
\begin{tabular}{rrrrrrrr}
\toprule
$B$ & block rows & encoded bits & QUBO terms & {\color{black}LU }$(2/3)B^3$ & rel. error & rel. residual & saturation / $\max\kappa_2(S_i)$ \\
\midrule
4  & 240 & 48  & 1,176 & $4.27\times10^{1}$ & $1.17\times10^{-3}$  & $3.08\times10^{-3}$  & $0$ / 2.96 \\
8  & 120 & 96  & 4,656 & $3.41\times10^{2}$ & $1.46\times10^{-3}$  & $8.75\times10^{-3}$  & $0$ / 4.04 \\
16 & 60  & 192 & 18,528& $2.73\times10^{3}$ & $2.35\times10^{-3}$  & $4.84\times10^{-2}$  & $0$ / 4.82 \\
32 & 30  & 384 & 73,920& $2.18\times10^{4}$ & $1.91\times10^{-3}$  & $1.20\times10^{-1}$  & $0$ / 5.29 \\
\bottomrule
\end{tabular}
\end{table*}

The block-dependent scaling rule yields zero observed clipping for every tested $B$, while the maximum Schur condition number rises only from 2.96 to 5.29. The relative solution error remains at or below $2.35\times10^{-3}$ through $B=32$. The relative residual is more sensitive, increasing to $1.20\times10^{-1}$ at $B=32$. This distinction is expected for elliptic operators: a visually and normwise accurate smooth field can have a noticeably larger residual because the differential operator amplifies high-frequency local perturbations.

The large-$B$ experiment tests fixed-point representability and downstream stability; optimization quality is examined separately at $B=8$.

% layout barrier relaxed

\subsection{Practical optimization of the 96-bit $B=8$ QUBOs}
\label{sec:results_phase4}

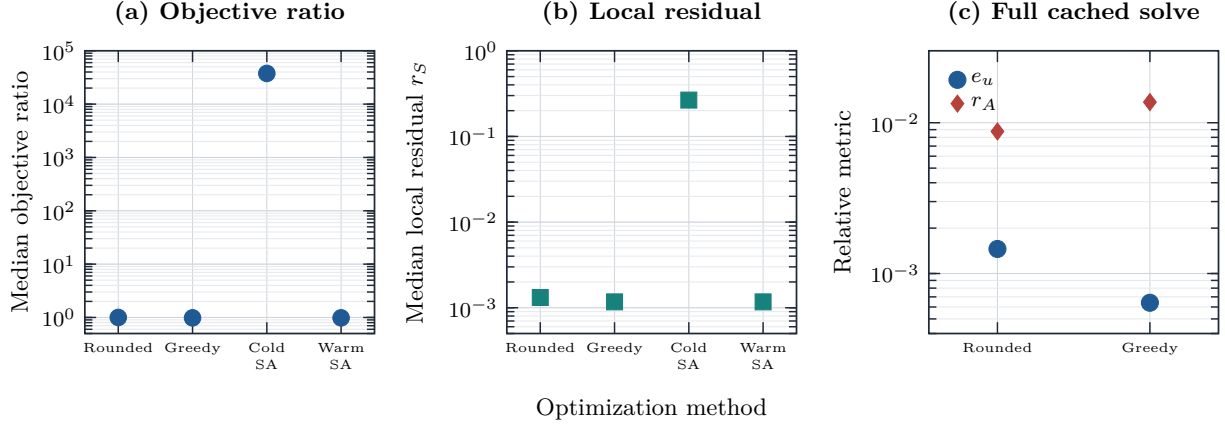
\begin{figure*}[t]
\centering
\begin{tikzpicture}
\begin{groupplot}[
    group style={group size=3 by 1, horizontal sep=1.75cm},
    qpde axis,
    width=0.295\textwidth,
    height=0.29\textwidth,
    x tick label style={rotate=0,anchor=north,align=center,font=\tiny},
    tick label style={font=\footnotesize},
]
\nextgroupplot[
 title={(a) Objective ratio},ymode=log,
 ylabel={\rev{Median objective ratio}},
 symbolic x coords={Rounded,Greedy,Cold SA,Warm SA},xtick=data,xticklabels={Rounded,Greedy,{\shortstack{Cold\\SA}},{\shortstack{Warm\\SA}}},enlarge x limits=0.15,
 ymin=5e-1,ymax=1e5]
\addplot[only marks,mark=*,mark size=3.1pt,color=qpdeblue] table[col sep=comma,x=short,y=median_objective_ratio_to_rounded]{phase4_plot.csv};
\nextgroupplot[
 title={(b) Local residual},ymode=log,
 ylabel={\rev{Median local residual $r_S$}},ylabel style={font=\small, at={(axis description cs:-0.27,0.5)}, anchor=south, inner sep=1pt},
 symbolic x coords={Rounded,Greedy,Cold SA,Warm SA},xtick=data,xticklabels={Rounded,Greedy,{\shortstack{Cold\\SA}},{\shortstack{Warm\\SA}}},enlarge x limits=0.15,
 ymin=5e-4,ymax=1]
\addplot[only marks,mark=square*,mark size=3.0pt,color=qpdeteal] table[col sep=comma,x=short,y=median_residual_norm]{phase4_plot.csv};
\nextgroupplot[
 title={(c) Full cached solve},ymode=log,
 ylabel={\rev{Relative metric}},ylabel style={font=\small, at={(axis description cs:-0.26,0.5)}, anchor=south, inner sep=1pt},
 symbolic x coords={Rounded,Greedy},xtick=data,enlarge x limits=0.45,
 ymin=4e-4,ymax=3e-2,
 legend style={at={(0.04,0.96)},anchor=north west,font=\footnotesize,legend columns=1}]
\addplot[only marks,mark=*,mark size=3.2pt,color=qpdeblue] coordinates {(Rounded,0.001455092)(Greedy,0.000640)};
\addlegendentry{\rev{$e_u$}}
\addplot[only marks,mark=diamond*,mark size=3.2pt,color=qpdered] coordinates {(Rounded,0.00875444)(Greedy,0.01370379)};
\addlegendentry{\rev{$r_A$}}
\end{groupplot}
\node[font=\small] at ($(group c1r1.south)!0.5!(group c3r1.south)+(0,-10mm)$) {\rev{Optimization method}};
\end{tikzpicture}
\caption{Practical optimization behavior for the $B=8$ Poisson line-block QUBOs. Panels (a) and (b) summarize the 32 representative 96-bit instances. In panel (a), the vertical axis is the median QUBO-objective ratio relative to direct fixed-point rounding, so a value below one indicates an improvement over the rounded baseline. In panel (b), the vertical axis is the median decoded local Schur-inverse residual $r_S$. Random cold-start simulated annealing lands far from the useful fixed-point region under the tested budget, whereas greedy and warm-start strategies preserve or modestly improve the rounded objective and maintain $10^{-3}$-level local residuals. Panel (c) uses all 960 inverse columns; its vertical axis shows the two end-to-end relative metrics, namely \rev{$e_u$, the final solution error, and $r_A$, the global PDE residual.} Greedy QUBO-space polishing substantially improves the final solution error, although the global PDE residual increases, demonstrating that a better collection of local objectives need not monotonically improve every global metric.}
\label{fig:phase4}
\end{figure*}

At $B=8$, each inverse-column QUBO contains 96 logical variables and 4,656 dense upper-triangular coefficients, so exhaustive enumeration would require $2^{96}$ states per column. The $B=8$ optimizer study tests whether practical QUBO-space heuristics recover states at least as useful as direct fixed-point rounding. Table~\ref{tab:phase4} and Figure~\ref{fig:phase4} show a sharp dependence on initialization.

The simulated-annealing budget is 140 sweeps per restart for cold starts and 120 for warm starts. Warm starts flip each rounded-state bit with probability $0.015$. Annealing uses geometric temperature decay to $10^{-8}$, followed by at most eight sweeps of best-improvement single-bit descent.

\begin{table*}[t]
\centering
\scriptsize
\caption{Representative 96-bit optimizer study on 32 inverse-column QUBOs from Schur blocks 1, 30, 60, and 120 of the $B=8$, $N=960$ Poisson problem. Objective ratios are computed relative to direct fixed-point rounding for the same column. The rounded state is a deterministic representability baseline, not an exact 96-bit QUBO optimum. Greedy descent and both annealing variants are followed by the same best-improvement single-bit polishing where applicable.}
\label{tab:phase4}
\begin{tabular}{lrrrrr}
\toprule
Method & median obj. ratio & mean obj. ratio & max obj. ratio & median $r_S$ & median column error \\
\midrule
Rounded fixed-point & 1.000 & 1.000 & 1.000 & $1.32\times10^{-3}$ & $3.26\times10^{-4}$ \\
Greedy from rounded & 0.983 & 0.900 & 1.000 & $1.17\times10^{-3}$ & $3.67\times10^{-4}$ \\
Cold SA + greedy & $3.76\times10^{4}$ & $1.33\times10^{5}$ & $1.50\times10^{6}$ & $2.66\times10^{-1}$ & $1.06\times10^{-1}$ \\
Warm SA + greedy & 0.983 & 0.881 & 1.000 & $1.17\times10^{-3}$ & $3.67\times10^{-4}$ \\
\bottomrule
\end{tabular}
\end{table*}

Cold random-start simulated annealing is ineffective under the tested budget: its median objective is $3.76\times10^4$ times the rounded objective and its median decoded local residual rises to $2.66\times10^{-1}$.
Across the 32 representative instances, the cold-start solutions have a
median Hamming distance of 47 bits from the rounded states; 31 retain at
least one sign- or integer-bit disagreement, and 31 undergo no further
flip during the final best-improvement single-bit polish. Since the sign
and integer bits contribute $-2\gamma_i$ and $+\gamma_i$, respectively,
recovering a low-residual inverse column from a random encoding can
require coordinated bit changes. The observed stagnation is consistent
with this local-search barrier under the tested budget. In contrast, greedy descent from the rounded state and warm-start annealing both remain close to the high-quality representable state. Their median objective ratio is $0.983$, and warm-start annealing reduces the mean objective ratio to $0.881$. The fixed-point inverse estimate provides an effective initialization for this dense 96-bit landscape.

The full-factor test strengthens this conclusion. Greedy single-bit polishing is applied to all 960 inverse columns. It finds an objective-decreasing modification for $48.3\%$ of the columns and lowers the mean local residual from $1.330\times10^{-3}$ to $1.285\times10^{-3}$. The final relative solution error improves from $1.455\times10^{-3}$ to $6.40\times10^{-4}$. The global PDE residual increases from $8.75\times10^{-3}$ to $1.37\times10^{-2}$, even as the final solution error decreases. Local least-squares improvement is therefore useful but not perfectly aligned with every global downstream metric, an observation that motivates residual-aware or globally informed polishing in future work.

% layout barrier relaxed

\subsection{Real QCI Dirac-3 hardware validation}
\label{sec:results_dirac}

The hardware experiment uses QCI Dirac-3 to optimize the 24-bit inverse-column QUBOs. All 34 unique submissions completed successfully. Jobs used two integer levels per variable, and no sum constraint. The upper-triangular QUBO coefficients were normalized by their maximum absolute magnitude before submission, with the constant term removed. Of the 34 jobs, 32 requested 20 samples at relaxation-schedule setting 1; two Poisson jobs (block 4, column 0 and block 6, column 1) requested 100 samples at setting 4. Because repeated Schur structures occur in the three PDE systems, these 34 jobs reconstruct 60 mapped inverse-column instances: 20 each for Heat, Poisson, and Klein--Gordon. The device returned the exact finite-bit reference bitstring for 31 of the 60 mapped instances and for 17 of the 34 unique submitted QUBOs. The total metered device time was 558~s. These timings document the hardware campaign but are not used to claim a speed advantage over classical PDE methods.

Table~\ref{tab:dirac} combines local hardware quality with the final reconstructed PDE accuracy. An exact hit occurs when the lowest-objective sample returned by Dirac-3 is identical to the reference bitstring obtained by {\color{black}exact reduced-grid minimization of the same 24-bit QUBO.} The exact-hit rate differs markedly across the three problems, from $20\%$ for Heat to $95\%$ for Klein--Gordon. All three reconstructed PDE solutions remain accurate. Bitstring recovery alone is therefore insufficient as an end-to-end metric: distinct bitstrings can decode to inverse columns with similarly small residuals.

\begin{table*}[t]
\centering
\scriptsize
\caption{QCI Dirac-3 hardware audit for the 60 mapped 24-bit inverse-column instances. Each QUBO is dense, with 300 upper-triangular coefficients. The normalized objective gap compares the best returned submitted objective with the exact finite-bit ground-state objective after the same normalization and constant removal. $\eta=r_S^{\rm D}/r_S^*$ is the local residual ratio relative to the exact finite-bit minimizer. The final two columns report the named end-to-end PDE reconstruction used in Figure~\ref{fig:dirac_reconstruction}.}
\label{tab:dirac}
\begin{tabular}{lrrrrrrr}
\toprule
PDE & exact hits & mean norm. gap & max norm. gap & mean $\eta$ & max $r_S^{\rm D}$ & PDE rel. error & PDE rel. residual \\
\midrule
Heat 1D & 4/20 & $2.40\times10^{-7}$ & $1.03\times10^{-6}$ & 2.11 & $1.04\times10^{-3}$ & $3.24\times10^{-4}$ & $3.33\times10^{-4}$ \\
Poisson 2D & 8/20 & $2.14\times10^{-6}$ & $1.03\times10^{-5}$ & 2.25 & $3.24\times10^{-3}$ & $1.97\times10^{-3}$ & $1.55\times10^{-2}$ \\
Klein--Gordon 1D & 19/20 & $1.94\times10^{-8}$ & $3.88\times10^{-7}$ & 1.07 & $6.87\times10^{-4}$ & $2.97\times10^{-4}$ & $2.97\times10^{-4}$ \\
\bottomrule
\end{tabular}
\end{table*}

\begin{figure*}[t]
\centering
\begin{tikzpicture}
\begin{groupplot}[
 group style={group size=3 by 1,horizontal sep=1.90cm},
 qpde axis,
 width=0.285\textwidth,height=0.265\textwidth,
 tick label style={font=\footnotesize},
 label style={font=\small},
 title style={font=\small\bfseries}]
\nextgroupplot[
 title={(a) Exact recovery},xlabel={\rev{PDE}},ylabel={\rev{Exact-recovery fraction}},ylabel style={font=\small, at={(axis description cs:-0.26,0.5)}, anchor=south, inner sep=1pt},
 symbolic x coords={Heat,Poisson,Klein--Gordon},xtick={Heat,Poisson,Klein--Gordon},xticklabels={Heat,Poisson,{\shortstack{Klein--\\Gordon}}},enlarge x limits=0.28,
 x tick label style={rotate=0,anchor=north,align=center,font=\scriptsize},ymin=0,ymax=1.05]
\addplot[ybar,bar width=11pt,bar shift=0pt,fill=qpdeblue!72,draw=qpdeblue] coordinates {(Heat,0.20)};
\addplot[ybar,bar width=11pt,bar shift=0pt,fill=qpdeteal!58,draw=qpdeteal] coordinates {(Poisson,0.40)};
\addplot[ybar,bar width=11pt,bar shift=0pt,fill=qpdeorange!48,draw=qpdeorange] coordinates {(Klein--Gordon,0.95)};
\nextgroupplot[
 title={(b) Residual ratio},xlabel={\rev{Mapped inverse-column instance}},ylabel={\rev{$\eta=r_S^{\rm D}/r_S^*$}},ylabel style={font=\small, at={(axis description cs:-0.24,0.5)}, anchor=south, inner sep=1pt},
 ymode=log,ymin=0.8,ymax=10,xmin=0,xmax=61]
\addplot[only marks,mark=*,mark size=1.7pt,color=qpdeblue] table[col sep=comma,x=idx,y=residual_ratio,restrict expr to domain={\thisrow{pde_code}}{1:1}]{dirac_scatter.csv};
\addplot[only marks,mark=square*,mark size=1.6pt,color=qpdeteal] table[col sep=comma,x=idx,y=residual_ratio,restrict expr to domain={\thisrow{pde_code}}{2:2}]{dirac_scatter.csv};
\addplot[only marks,mark=triangle*,mark size=1.9pt,color=qpdeorange] table[col sep=comma,x=idx,y=residual_ratio,restrict expr to domain={\thisrow{pde_code}}{3:3}]{dirac_scatter.csv};
\addplot[qpde reference] coordinates {(1,1)(60,1)};
\nextgroupplot[
 title={(c) Gap vs. residual},xlabel={\rev{Normalized objective gap $g_j$}},ylabel={\rev{$r_S^{\rm D}$}},ylabel style={font=\small, at={(axis description cs:-0.28,0.5)}, anchor=south, inner sep=1pt},
 xmode=log,ymode=log,xmin=5e-10,xmax=5e-5,ymin=7e-5,ymax=8e-3]
\addplot[only marks,mark=*,mark size=1.7pt,color=qpdeblue] table[col sep=comma,x expr={max(\thisrow{normalized_objective_gap},1e-9)},y=local_inverse_residual_dirac,restrict expr to domain={\thisrow{pde_code}}{1:1}]{dirac_scatter.csv};
\addplot[only marks,mark=square*,mark size=1.6pt,color=qpdeteal] table[col sep=comma,x expr={max(\thisrow{normalized_objective_gap},1e-9)},y=local_inverse_residual_dirac,restrict expr to domain={\thisrow{pde_code}}{2:2}]{dirac_scatter.csv};
\addplot[only marks,mark=triangle*,mark size=1.9pt,color=qpdeorange] table[col sep=comma,x expr={max(\thisrow{normalized_objective_gap},1e-9)},y=local_inverse_residual_dirac,restrict expr to domain={\thisrow{pde_code}}{3:3}]{dirac_scatter.csv};
\end{groupplot}
\end{tikzpicture}
\caption{Audit of the 60 mapped inverse-column instances reconstructed from 34 unique QCI Dirac-3 submissions. Panel (a) reports exact finite-bit ground-state recovery, with the vertical axis giving the fraction of mapped instances for which the hardware returned the exact finite-bit minimizer. \rev{Its Heat, Poisson, and Klein--Gordon bar colors match the corresponding marker colors in panels (b) and (c).} Panel (b) compares each hardware-decoded local residual with the residual of the exact finite-bit minimizer; its vertical axis is the residual ratio $\eta=r_S^{\rm D}/r_S^*$, where $r_S^{\rm D}$ is the residual obtained from the Dirac-3 decoded solution and $r_S^*$ is the residual of the exact finite-bit minimizer. Thus, $\eta=1$ corresponds to equal residual, \rev{as indicated by the horizontal dashed reference line}. Panel (c) relates the normalized submitted-objective gap $g_j$ on the horizontal axis to the absolute hardware-decoded local residual $r_S^{\rm D}$ on the vertical axis; exact hits are displayed at the plotting floor $10^{-9}$ on the logarithmic horizontal axis. The figure shows that exact bitstring recovery is neither uniform across PDEs nor necessary for a numerically useful cached inverse column.}
\label{fig:dirac_audit}
\end{figure*}
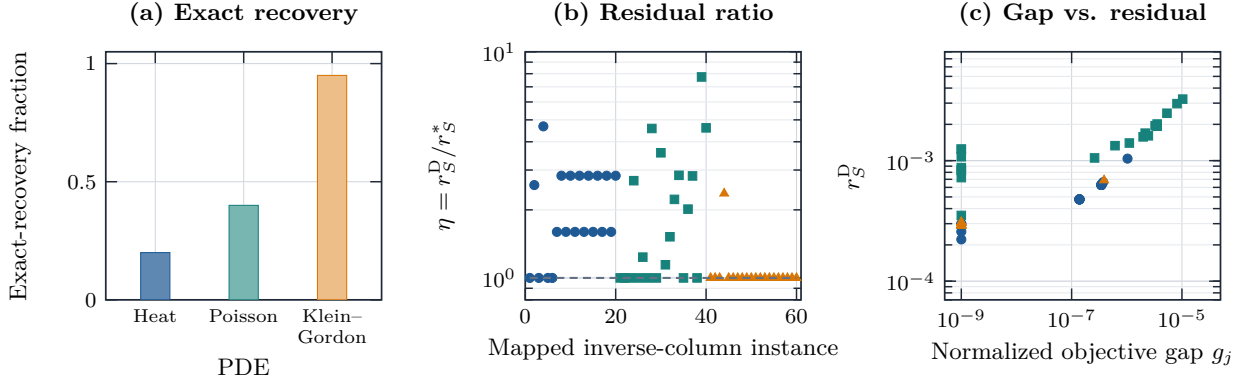

Figure~\ref{fig:dirac_audit} makes this distinction visible. Panel (a) gives the exact finite-bit recovery fraction. Panel (b) shows the local residual ratio for all mapped instances: non-exact Heat and Poisson samples often increase the local residual by a factor of a few, but the absolute residuals remain of order $10^{-3}$ or smaller. Panel (c) compares normalized objective gap with the decoded local residual. Exact hits lie at zero gap, while the non-exact samples occupy a narrow low-gap region.

% layout barrier relaxed

\begin{figure*}[t]
\centering
\begin{tikzpicture}
\begin{groupplot}[
 group style={group size=2 by 2,horizontal sep=1.85cm,vertical sep=1.65cm},
 qpde axis,
 width=0.405\textwidth, enlarge x limits=0.04,
 enlarge y limits=0.06,height=0.25\textwidth,
 tick label style={font=\small},
 label style={font=\small},
 title style={font=\small\bfseries},
 legend style={font=\footnotesize},
]
\nextgroupplot[title={(a) Heat 1D},xlabel={\rev{Position $x$}},ylabel={\rev{$u(x)$}},legend pos=south east]
\addplot[qpde orange,solid,mark=o,mark size=1.5pt] table[col sep=comma,x expr={\thisrow{index}/21},y=qci_quantum]{dirac_heat.csv};\addlegendentry{Dirac-3 factors}
\addplot[black!65,very thick,densely dashed,mark=none] table[col sep=comma,x expr={\thisrow{index}/21},y=dense_ground_truth]{dirac_heat.csv};\addlegendentry{dense reference}
\nextgroupplot[title={(b) Poisson 2D slices},xlabel={\rev{Position $x$}},ylabel={\rev{$u(x_i,y_j)$}}]
\addplot[qpde orange,solid,mark=o,mark size=1.4pt] table[col sep=comma,x expr={\thisrow{x_index}/11},y=hardware,restrict expr to domain={\thisrow{slice}}{1:1}]{dirac_poisson_slices.csv};
\addplot[black!65,very thick,densely dashed,mark=none] table[col sep=comma,x expr={\thisrow{x_index}/11},y=dense,restrict expr to domain={\thisrow{slice}}{1:1}]{dirac_poisson_slices.csv};
\addplot[qpde red,dashdotted,mark=square*,mark size=1.3pt] table[col sep=comma,x expr={\thisrow{x_index}/11},y=hardware,restrict expr to domain={\thisrow{slice}}{2:2}]{dirac_poisson_slices.csv};
\addplot[black!50,very thick,densely dotted,mark=none] table[col sep=comma,x expr={\thisrow{x_index}/11},y=dense,restrict expr to domain={\thisrow{slice}}{2:2}]{dirac_poisson_slices.csv};
\nextgroupplot[title={(c) Klein--Gordon: $u$},xlabel={\rev{Position $x$}},ylabel={\rev{Displacement $u(x)$}}]
\addplot[qpde orange,solid,mark=o,mark size=1.5pt] table[col sep=comma,x expr={\thisrow{space_index}/11},y=u_qci_quantum]{dirac_kg.csv};
\addplot[black!65,very thick,densely dashed,mark=none] table[col sep=comma,x expr={\thisrow{space_index}/11},y=u_dense_ground_truth]{dirac_kg.csv};
\nextgroupplot[title={(d) Klein--Gordon: $v$},xlabel={\rev{Position $x$}},ylabel={\rev{Velocity $v(x)$}}]
\addplot[qpde orange,solid,mark=o,mark size=1.5pt] table[col sep=comma,x expr={\thisrow{space_index}/11},y=v_qci_quantum]{dirac_kg.csv};
\addplot[black!65,very thick,densely dashed,mark=none] table[col sep=comma,x expr={\thisrow{space_index}/11},y=v_dense_ground_truth]{dirac_kg.csv};
\end{groupplot}
\end{tikzpicture}
\caption{End-to-end PDE reconstructions using cached inverse factors decoded from real QCI Dirac-3 samples. Panel (a) compares the Heat 1D field with the dense reference. Panel (b) displays the $y$-line slice of the $10\times2$ Poisson interior grid; \rev{the plotted curve represents the physical coordinate $x_i=i/11$, with the slice at $y_1=1/3$.} Solid/dashed styles distinguish hardware and reference traces. Panels (c) and (d) show the two components of the interleaved Klein--Gordon state on \rev{$x_i=i/11$}. The corresponding relative solution errors are $3.24\times10^{-4}$ (Heat), $1.97\times10^{-3}$ (Poisson), and $2.97\times10^{-4}$ (Klein--Gordon).}
\label{fig:dirac_reconstruction}
\end{figure*}
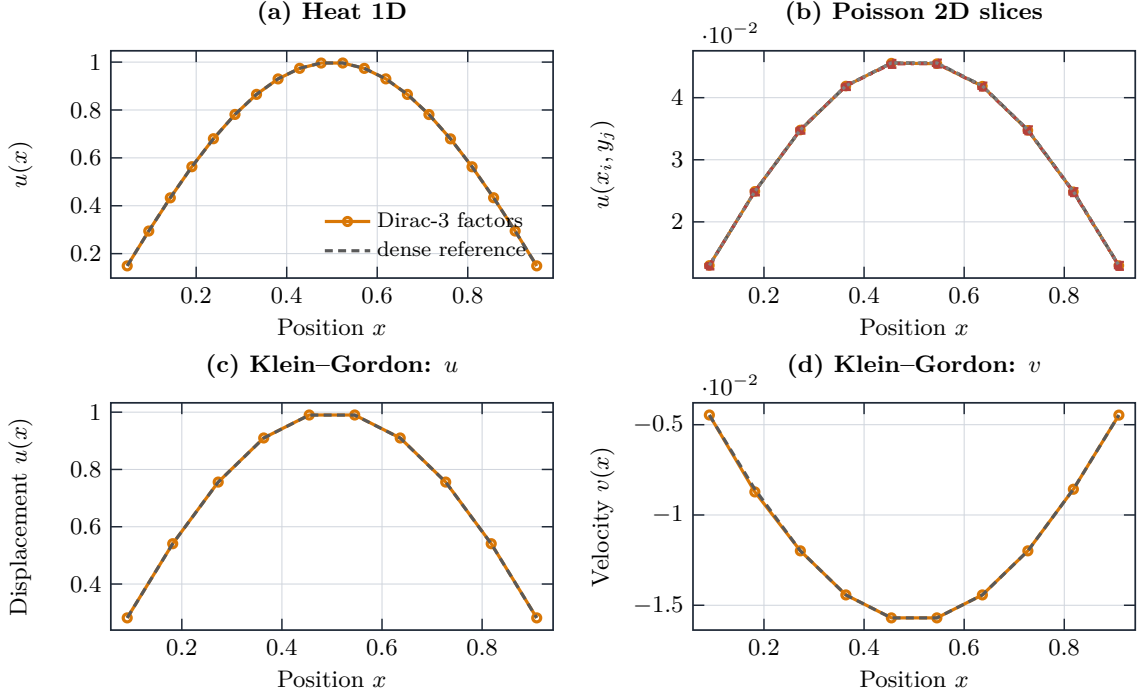

The end-to-end reconstructions are shown in Figure~\ref{fig:dirac_reconstruction}. For Heat, the relative error is $3.24\times10^{-4}$. For Poisson, the two line slices of the $10\times2$ interior grid are both reproduced closely; the relative field error is $1.97\times10^{-3}$ even though the relative residual is larger, $1.55\times10^{-2}$. For Klein--Gordon, both components of the interleaved state follow the dense reference with a combined relative error of $2.97\times10^{-4}$. The hardware experiment executes the complete interface. PDE-derived QUBOs are submitted to the device, decoded into local inverse columns, assembled into cached factors, and used in the classical PDE solve.

\section{Discussion and limitations}
\label{sec:discussion}

The experiments establish that a fixed-operator finite-difference solve can be reorganized so that QUBO optimization constructs reusable local inverse factors rather than each global PDE solution. The mathematical consequence is that the largest logical QUBO has dimension $B(1+M+K)$, controlled by local block size and precision, while the number of blocks determines how many such problems must be solved and cached. The numerical consequence is an offline--online separation.

The validation proceeds in four stages. The exact 24-bit studies establish the finite-point theory without optimizer ambiguity and recover the predicted $O(2^{-K})$ inverse/solution scaling and $O(4^{-K})$ squared-residual scaling. The five-PDE tests show that the same cached-factor construction works for parabolic, nonlinear IMEX, hyperbolic, elliptic, and oscillatory finite-difference systems using predetermined multiscale correction grids. The $B$ sweep then shows that the representation remains usable for substantially larger local blocks when block-dependent scaling prevents clipping. The 96-bit experiment reveals the optimization difficulty that representability alone cannot address: random cold starts are poor, whereas a structure-aware rounded or warm-start state provides a useful search region. Finally, the Dirac-3 campaign demonstrates that the local QUBO interface can be executed on QCI Dirac-3 hardware and that the returned states can be propagated all the way back to accurate PDE fields.

The present study establishes the architecture and its numerical accuracy. Performance relative to sparse LU, block Thomas, Krylov, and multigrid solvers remains a separate benchmark question. The contribution is architectural. It identifies an operator-dependent preprocessing kernel that can be exposed to a QUBO optimizer and distributed across repeated solves. Whether such an interface becomes computationally advantageous depends on future optimizer and hardware scaling as well as on the number of right-hand sides over which preprocessing can be reused.

Reducing logical variable count leaves the dense interaction cost intact. For the line-block Poisson experiment, the encoded dimension grows linearly with $B$, but the least-squares QUBO is fully dense and its coefficient count grows quadratically. On hardware with sparse native connectivity this may introduce embedding, routing, chain, or precision overhead. The present Dirac-3 interface accepts the dense polynomial form directly, but the same logical QUBO can have very different physical costs on other architectures. Logical bits, QUBO coefficients, physical variables, and elapsed device time should therefore not be treated as interchangeable resource measures.

The method requires a block-tridiagonal or compatible block-banded ordering. Many finite-difference grids admit such orderings naturally, but irregular meshes and general sparse operators may require nested dissection, domain decomposition, cyclic reduction, low-rank interface corrections, or a hierarchical Schur construction. The scaling parameter also remains problem dependent. The per-block diagonal-dominance rule used in the large-$B$ experiment eliminates clipping for the tested Poisson operators, but nonsymmetric, indefinite, or weakly diagonally dominant Schur blocks may require a more general adaptive range estimator.

The hardware campaign uses 60 mapped 24-bit instances reconstructed from 34 unique Dirac-3 submissions. The 24-bit problems are small enough to retain exact classical references. Larger blocks, higher fixed-point precision, larger two- and three-dimensional grids, multiple hardware platforms, and statistically repeated device runs are required before drawing conclusions about hardware scaling. The current hardware data support the end-to-end interface demonstrated here.

Two extensions follow directly from the present results. On the numerical side, adaptive correction-scale selection, block-ordering strategies, and integration with Krylov or multigrid outer solvers can make the cached factors useful as approximate inverses or preconditioners in larger applications. On the optimization side, the same inverse-column benchmark family provides a controlled way to compare classical heuristics, quantum-inspired solvers, annealers, and other quantum optimization hardware at increasing $B$ while preserving a clear downstream PDE accuracy test.

\section{Conclusions}
\label{sec:conclusion}

We introduced a reusable block-Schur formulation for fixed-operator finite-difference PDE systems in which QUBO optimization is confined to the offline construction of local inverse factors. Each Schur-inverse column is represented on a signed fixed-point grid, obtained from a least-squares QUBO, decoded, and cached. Subsequent right-hand sides are solved classically with the stored factors. This construction separates the logical QUBO dimension from the global PDE dimension: at fixed precision the largest logical QUBO contains $B(1+M+K)$ variables, although the number of local problems and their coupling density remain important resource costs.
The results establish reusable local QUBO kernels as an interface between specialized discrete optimization and conventional numerical PDE workflows. {\color{black}Assessing scalability and any quantum computational advantage over established classical PDE solvers requires larger-scale experiments, and comparative performance benchmarks remain a direction for future work.}

\FloatBarrier
\section*{Data and code availability}
\label{sec:data_availability}
The code, numerical data, processed datasets, hardware audit records, and scripts used to reproduce the analyses and figures are publicly available at \url{https://github.com/Ar2000Ash/QPDE}. The repository provides environment requirements, experiment drivers, figure-generation scripts, and numerical verification tests. The QCI results can be verified offline from recorded device outputs; physical hardware execution requires access to the corresponding optimization platform.

\section*{Declaration of competing interest}

The authors declare that they have no known competing interests.

\section*{Funding}

This research received no external funding.

\section*{AI disclosure}

OpenAI ChatGPT (GPT-5.6 Sol) was used to assist with manuscript organization, presentation, and language editing. The authors take full responsibility for the content, accuracy, interpretation, and conclusions of the work.

\bibliographystyle{elsarticle-num}
\bibliography{refs}

\end{document}

%% file: qpde_figure_style.tex
\usepackage{tikz}
\usepackage{pgfplots}
\usepackage{pgfplotstable}
\usepackage{siunitx}
\usetikzlibrary{arrows.meta,positioning,fit,calc,backgrounds,matrix,shapes.geometric}
\usepgfplotslibrary{groupplots}
\pgfplotsset{compat=1.18}

\definecolor{qpdeblue}{HTML}{1F5A94}
\definecolor{qpdeteal}{HTML}{17827B}
\definecolor{qpdeorange}{HTML}{D97706}
\definecolor{qpdered}{HTML}{B6403A}
\definecolor{qpdegray}{HTML}{667085}
\definecolor{qpdedark}{HTML}{1F2937}
\definecolor{qpdeline}{HTML}{D0D5DD}
\definecolor{qpdelight}{HTML}{F5F7FA}

\pgfplotsset{
  qpde axis/.style={
    width=\linewidth,
    height=0.64\linewidth,
    axis line style={draw=qpdedark, line width=0.65pt},
    tick style={draw=qpdedark, line width=0.55pt},
    tick label style={font=\small},
    label style={font=\small},
    xlabel style={font=\small, yshift=1pt},
    ylabel style={font=\small, at={(axis description cs:-0.17,0.5)}, anchor=south, inner sep=1pt},
    title style={font=\small\bfseries},
    legend style={
      font=\footnotesize,
      draw=none,
      fill=none,
      cells={anchor=west},
      /tikz/every even column/.append style={column sep=0.55em}
    },
    grid=both,
    major grid style={draw=qpdeline, line width=0.35pt},
    minor grid style={draw=qpdeline!55, line width=0.25pt},
    enlargelimits=false,
    clip marker paths=true,
  },
  qpde blue/.style={color=qpdeblue, mark=*, very thick, mark size=1.8pt},
  qpde teal/.style={color=qpdeteal, mark=square*, very thick, mark size=1.7pt},
  qpde orange/.style={color=qpdeorange, mark=triangle*, very thick, mark size=2.0pt},
  qpde red/.style={color=qpdered, mark=diamond*, very thick, mark size=1.8pt},
  qpde reference/.style={color=qpdegray, densely dashed, thick, mark=none},
}

\tikzset{
  qpdebox/.style={
    rounded corners=2.2pt,
    draw=qpdedark!82,
    line width=0.7pt,
    fill=white,
    minimum height=10mm,
    inner xsep=3mm,
    inner ysep=2mm,
    align=center,
    font=\small
  },
  qpdeaccent/.style={qpdebox, draw=qpdeblue, fill=qpdeblue!5},
  qpdecache/.style={qpdebox, draw=qpdeteal, fill=qpdeteal!6},
  qpdearrow/.style={-{Latex[length=2.2mm,width=1.5mm]}, line width=0.8pt, draw=qpdedark!85},
  qpdedashedarrow/.style={qpdearrow, dashed, draw=qpdegray},
  qpdepanel/.style={font=\small\bfseries, text=qpdedark},
}